\documentclass[12pt]{amsart}
\usepackage{amsfonts}
\usepackage{amsthm}
\usepackage{amsmath}    
\usepackage{graphicx}   
\usepackage{verbatim}   
\usepackage{color}      
\usepackage{subfigure}  
\usepackage{epsfig}
\usepackage{amssymb}
\usepackage{mathrsfs}
\usepackage{tikz}

\theoremstyle{definition}

\newtheorem{rem}{Remark}[section]

\def\vn{{\boldsymbol v}}

\def\q{{\boldsymbol q}}
\def\n{{\boldsymbol n}}

\def\x{{\boldsymbol x}}
\def\p{{\boldsymbol p}}

\def\s{{\boldsymbol s}}

\def\psin{{\boldsymbol \psi }}
\def\etan{{\boldsymbol \eta }}
\def\rhon{{\boldsymbol  \rho }}
\def\zetan{{\boldsymbol  \zeta }}
\def\nun{{\boldsymbol  \nu }}

\newcommand{\abs}[1]{\left\vert#1\right\vert}

\def\0{{\boldsymbol 0}}

\def\n{{\boldsymbol n}}
\def\n{{\boldsymbol n}}
\def\p{{\boldsymbol p}}
\def\s{{\boldsymbol s}}

\def\x{{\boldsymbol x}}
\def\vn{{\boldsymbol v}}
\def\r{{\boldsymbol r}}
\def\etan{{\boldsymbol \eta}}
\def\thetan{{\boldsymbol \theta}}
\def\psin{{\boldsymbol \psi}}
\def\rhon{{\boldsymbol \rho}}
\def\zetan{{\boldsymbol \zeta}}

\def\beq{\begin{equation}}
\def\eeq{\end{equation}}

\makeatletter      
\@addtoreset{equation}{section}
\makeatother       
\makeatletter      
\@addtoreset{table}{section}
\makeatother       
\makeatletter      
\@addtoreset{figure}{section}
\makeatother       
\allowdisplaybreaks

\def\signFF{\bigskip\bigskip\hspace{80mm}
\vbox{{\sc Francis Filbet\par\vspace{3mm}
Universit\'e de Toulouse III \& IUF \par
UMR5219, Institut de Math\'ematiques de Toulouse,\par
118, route de Narbonne\par
F-31062 Toulouse cedex,  FRANCE
\par\vspace{3mm}e-mail:} francis.filbet@math.univ-toulouse.fr }}

\def\signRG{\bigskip\bigskip\hspace{80mm}
\vbox{{\sc Ruihan Guo \par\vspace{3mm}
Institut Camille Jordan,\par
Universit\'{e} Claude Bernard Lyon
I, \par
43 boulevard 11 novembre 1918, \par
69622 Villeurbanne cedex, France.
\par\vspace{3mm}e-mail:}  rguo@math.univ-lyon1.fr. }}

\begin{document}

\title[Local Discontinuous Galerkin level set method for Willmore flow]{A $p$-adaptive local discontinuous Galerkin level set method for Willmore flow}%

\author{Ruihan Guo and Francis Filbet}

\date{\today}
\maketitle

\begin{abstract}
The level set method is often used to capture interface behavior in two or three dimensions.  In this paper, we present a combination of local discontinuous Galerkin (LDG) method and level set method for simulating Willmore flow. The LDG scheme is energy stable and mass conservative, which are good properties comparing with other numerical methods. In addition, to enhance the efficiency of the proposed LDG scheme and level set method, we employ a $p$-adaptive local discontinuous Galerkin technique, which applies high order polynomial approximations around the zero level set and low order ones away from the zero level set. A major advantage of the level set method is that the topological changes are well defined and easily performed. In particular, given the stiffness of Willmore flow, a high order semi-implicit Runge-Kutta method is employed for time discretization, which allows larger time step. The equations at the implicit time level are linear, we demonstrate an efficient and practical multigrid solver to solve the equations.  Numerical examples are given to illustrate the combination of the LDG scheme and level set method provides an efficient and practical approach when simulating the Willmore flow.

\smallskip

\bigskip

\end{abstract}

\tableofcontents

\section{Introduction}
\label{se:Intro}
The level set method was first introduced by Osher and Sethian \cite{osher} in 1987, and it has had large impact on computational methods for interface motion. In the past years, the level set method attracted a lot of attention and was used for a wide variety of problems, including fluid dynamics, material sciences, computer vision and imaging processing. For an extensive exposition of the level set method, we refer readers to the review paper \cite{osher_review} and the book \cite{osher_book}.

In the level  set method, the evolving family of interface $\Gamma_t$, $t\geq 0$ is represented implicitly as the zero level set of a continuous function which we will denote as $\phi: \Omega \times [0,T]\rightarrow \mathbb{R}$, where $\Omega \subset \mathbb{R}^2$ is a simply connected domain containing the whole family of evolving interfaces $\Gamma_t$, $t \in [0,T]$. Therefore, the interface $\Gamma_t$ can be written as
\begin{equation*}
\Gamma_t\,=\,\left\{\x\,\in\,\mathbb{R}^2,\quad \phi(\x,t)=0\,\right\}.
\end{equation*}
There are certain advantages associated with this implicit Eulerian type of representation.
\begin{itemize}
\item The level set formulation is unchanged in higher dimensions, that is, it is easy to extend to 3D.
\item Topological changes such as merging and breaking in the evolving interface $\Gamma_t$ are handled naturally.
\end{itemize}

The level set formulation for Willmore flow will typically result in a fourth order nonlinear partial differential equation (PDE), {\it e.g.} see \cite{Droske} for details.
Bene\v{s} {\it et al.} developed a finite volume method for spatial discretization in \cite{Bene}. When the evolved surface is a graph, the Willmore flow converts to Willmore flow of graphs, and numerical methods for Willmore flow of graphs attracted a lot of attention in the last decades. A semi-implicit numerical scheme for the Willmore flow of graphs with continuous finite element discretization and the convergence analysis has been provided by Deckelnick and Dziuk in \cite{Deck}. Finite difference discretization for the Willmore flow of graphs has been discussed in \cite{Ober}. Xu and Shu \cite{Xu_will} developed a local discontinuous Galerkin method for Willmore flow of graphs, and time discretization was by the forward Euler method with a suitably small time step for stability.

The discontinuous Galerkin (DG) method is a class of finite element methods using completely discontinuous basis functions, which are usually chosen as piecewise polynomials. It was first introduced in 1973 by Reed and Hill \cite{Reed} for solving the linear time dependent neutron transport equation. Later, the DG method was developed for solving the nonlinear hyperbolic conservation laws by Cockburn {\it et al.} in a series of papers \cite{Cockburn3, Cockburn4, Cockburn5, Cockburn6}.

These discontinuous Galerkin methods were generalized to solve the
PDEs containing higher than first order derivatives, {\it i.e.} the
local discontinuous Galerkin (LDG) method. The first LDG method was
designed to solve a convection diffusion equation (with second
derivatives) by Cockburn and Shu \cite{Cockburn_Shu}. The idea of LDG
methods is to suitably rewrite a higher order PDE into a first order
system, then apply the DG method to the system. A key ingredient for
the success of such methods is the correct design of interface
numerical fluxes, which should be designed to guarantee stability. For
a detailed description about the LDG methods for high order time
dependent PDEs, we refer the readers to \cite{xu_review}. The DG and
LDG methods have several attractive properties, for example:
\begin{enumerate}
\item The order of accuracy can be locally determined in each cell,
  thus allowing for efficient $p$ adaptivity.
\item These methods can be used on arbitrary triangulations, even
  those with hanging nodes, thus allowing for efficient $h$
  adaptivity.
\item These methods have excellent parallel efficiency.
\end{enumerate}

In this paper, we present the level set method for computing evolution of the Willmore flow of embedded surfaces, and one of the major  advantages of level set method is their ability to easily handle topological changes. Then we develop a high order local discontinuous Galerkin method for solving level set formulation of the Willmore flow. The basis functions of LDG method can be completely discontinuous, which leads that the order of accuracy can be locally determined in each cell, {\it i.e.} $p$-adaptivity. To improve the efficiency of the proposed methods, a $p$-adaptive technique is applied here to  capture the evolution of the interface.

Another main difficulty when solving the Willmore flow numerically is
that it is stiff, and explicit time discretization methods suffer from
extremely small time step restriction for stability, of the form
$\Delta t \leq C \Delta x^4$, which is not efficient, especially for
long time simulation.  It would therefore be desirable to develop
implicit or semi-implicit time marching techniques to alleviate this
problem. In \cite{Xia_time}, three efficient time discretization
techniques were explored, these are the semi-implicit spectral
deferred correction (SDC) method, the additive Runge-Kutta (ARK)
method and the exponential time differencing (ETD) method. These
methods are mainly based on a splitting of stiff and non stiff
differential operators, which is particularly efficient for
semi-linear PDEs, but for quasi-linear equations, like  for Willmore
flows, fully nonlinear solver would be required. Thus we will employ
here a high order semi-implicit Runge-Kutta method \cite{Filbet, RFX}.

The organization of the paper is as follows. In Section \ref{se:level}, we review the level set formulation for Willmore flow and get a highly nonlinear fourth order PDE. In Section \ref{se:ldg}, we present the formulation of the local discontinuous Galerkin method and the $p$-adaptive technique for solving the Willmore flow. Section \ref{se:time} contains a simple description of the high order semi-implicit Runge-Kutta method. Numerical experiments are presented in Section \ref{se:numerical}, testing the performance of $p$-adaptive local discontinuous Galerkin level set method for Willmore flow. Finally, we give concluding remarks in Section \ref{se:con}.

\section{Level set formulation for Willmore flow}
\label{se:level}
In the level  set method, the evolving family of interface $\Gamma_t$, $t\geq 0$ is represented implicitly as the zero level set of a continuous function  $\phi$.
Therefore, the interface $\Gamma_t$ can be written as
\begin{equation*}
\Gamma_t\,=\,\left\{\x\in\mathbb{R}^2, \quad\phi(\x,t)=0\right\}.
\end{equation*}
For the Willmore flow, the level set function
$$
\frac{\partial \phi}{\partial t}\,+\,{\bf V}\,|\nabla \phi|\, \,=\,\,0,
$$
with ${\bf V}$ being the speed of propagation of the level set function $\phi$, turns into the following equation
\begin{equation}
\frac{\partial \phi}{\partial t}\,+\,\left(\Delta_{\Gamma_t} H(t)+H(t)\,\left(\|S(t)\|_{L^2}^2\,-\,\frac{1}{2}\,H(t)^2\right)\right)\,|\nabla \phi|\,=\,0,
\end{equation}
where $S(t)$ denotes the shape operator on $\Gamma_t$, $H(t)$ is the mean curvature on $\Gamma_t$ and is given as
\begin{equation*}
H=\nabla \cdot \left(\frac{\nabla \phi}{|\nabla \phi|}\right).
\end{equation*}
In order to obtain the Willmore flow, it is necessary to reformulate the equation in terms of $\phi(t)$, $H(t)$ and their derivatives.

For the sake of simplicity, we assume that $|\nabla \phi| \neq 0$ on $\Gamma_t$. Then the unit normal vector to a level curve of $\phi$ can be defined as
\begin{equation}
\n=\frac{\nabla \phi}{|\nabla \phi|}.
\end{equation}
Let us denote the following auxiliary functions $Q\,=\,|\nabla \phi|$
and $\omega\,=\,Q\,H$. Then, by the derivation in \cite{Droske}, a level set formulation of Willmore flow is given as follows:
\begin{equation}
\label{willmore}
\left\{
\begin{array}{l}
\displaystyle{\frac{\partial \phi}{\partial t}\,+\, Q \,\nabla \cdot \left(\mathbb{E}\,
  \nabla \omega\,-\,\frac{1}{2}\,\frac{\omega^2}{Q^3}\,\nabla \phi\right)\,=\,0,}
\\
\;
\\
\displaystyle{\omega =Q\,\nabla \cdot\left(\frac{\nabla \phi}{|\nabla \phi|}\right).}
\end{array}\right.
\end{equation}
where the matrix $\mathbb{E}$ is given by
\begin{equation*}
\mathbb{E}\,=\,\frac{1}{Q}\,\left(\mathbb{I}-\frac{\nabla \phi}{Q}\otimes \frac{\nabla \phi}{Q}\right),
\end{equation*}
which  is a projection into a tangential space of the curve representing the zero level set of $\phi$.

 The Willmore flow \eqref{willmore} is subjected to the initial condition
 \begin{equation*}
 \phi(\x,0)\,=\,\phi^0(\x), \quad \x \,\in\, \Omega,
 \end{equation*}
 where the initial function $\phi^0(\x)$ is a signed distance
 function, {\it i.e.} $\phi^0(\x)=\text{dist}(x,\Gamma_0)$, for any $\x \in \Omega$. The initial condition satisfies the following property
 \begin{equation}
 |\nabla \phi^0(\x)|\,=\,1.
 \end{equation}

 One of the major advantages of level set methods is their ability to easily handle topological changes. However, we have found this not to be the case for Willmore flow because of the lack of a maximum principle. Indeed, two surfaces both undergoing an evolution by Willmore flow may intersect in finite time. Therefore, the level set formulation for Willmore flow will lead to singularities in general, {\it i.e.} $\nabla \phi=0$ in finite time. To handle these difficulties, it is necessary to introduce a suitable regularization technique. For a fixed regularization parameter $0<\varepsilon \ll 1$, we define
 \begin{equation}
 Q_{\varepsilon}\,=\,\sqrt{\varepsilon+|\nabla \phi|^2}
 \end{equation}
and replace all occurrences of $Q$ in \eqref{willmore} by $Q_{\varepsilon}$.

 \section{Numerical method}
 \label{se:ldg}
With the level set formulation for Willmore flow, we need to discrete equation \eqref{willmore} to capture the movement of the interface (zero level set). It would  therefore be desirable to develop high order numerical scheme in both space and time to obtain high resolution. In this section, we will consider the local discontinuous Galerkin (LDG) method for the Willmore flow \eqref{willmore}.
\subsection{Notation}
We consider a subdivision $\mathcal{T}_h$ of $\Omega$ with shape-regular elements $K$. Let $\mathcal{E}_h$ denote the union of the boundary faces of elements $K \in \mathcal{T}_h$,
{\it i.e.} $\mathcal{E}_h=\cup_{K\in \mathcal{T}_h} \partial K$, and $\mathcal{E}_0=\mathcal{E}_h \setminus  \partial \Omega$.

In order to describe the flux
functions, we need to introduce some notations. Let $e$ be an interior face
shared by the ``left'' and ``right'' elements $K_L$ and $K_R$ and define the normal vectors $\nun_L$ and
$\nun_R$ on $e$ pointing exterior to $K_L$ and $K_R$, respectively. For
our purpose ``left'' and ``right'' can be uniquely defined for each
face according to any fixed rule. For example, we choose $\nun_0$ as a constant vector. The left element $K_L$ to the face $e$ requires that $\nun_L \cdot \nun_0<0$, and the right one $K_R$ requires $\nun_L \cdot \nun_0 \geq 0$.
If $\psi$ is a function on  $K_L$ and $K_R$, but possibly
discontinuous across $e$, let $\psi_L$ denote $(\psi|_{K_L})|_{e}$
and $\psi_R$ denote $(\psi|_{K_R})|_{e}$, the left and right trace,
respectively.

Let $\mathcal{P}^k(K)$ be the
space of polynomials of degree at most $k \geq 0$ on $K \in
\mathcal{T}_h$. The finite element spaces are denoted by
$$
V_{h}\,=\,\Bigl\{\varphi\in L^2(\Omega),\quad \varphi\vert_K\in
\mathcal{P}^k(K),\quad \forall \,K\in \mathcal{T}_{h}\Bigr\}
$$
and
$$
\Sigma_h^d \,=\,  \Bigl\{\Phi=(\phi_1,\cdots,\phi_d)^T \in
(L^2(\Omega))^d, \quad \phi_l\vert_{ K }\in \mathcal{P}^k(K), \,
\quad \forall l=1\cdots d,\,\forall \,K \in \mathcal{T}_{h}\Bigr\}.
$$
Note that functions in $V_{h}$, $\Sigma_h^d$  are allowed to
be completely discontinuous across element interfaces.
\subsection{Local discontinuous Galerkin method}
To develop the local discontinuous Galerkin scheme, we first rewrite the Willmore flow \eqref{willmore} as a first order system:
\begin{equation}
\label{eq:willmore-sys}
\left\{
\begin{array}{l}
\displaystyle{\frac{1}{Q_\varepsilon}\,\frac{\partial \phi}{\partial t} \,+\,\nabla
\cdot (\s-\vn)\,=\,0 ,}
\\
\,
\\
\s =\mathbb{E}\, \p \quad{\rm and}\quad\displaystyle{\vn =\frac{1}{2} \frac{\omega^2}{Q_\varepsilon^3} \r,}
\end{array}\right.
\end{equation}
with $\p =\nabla \omega$, where $\omega$ is given by the algebraic
equation $\omega =Q_\varepsilon H$, and $H$ is solution to $H =\nabla \cdot
\q$, where $\q$ satisfies $Q_\varepsilon\,\q =\r$ and finally $\r=
\nabla \phi$. In equation \eqref{eq:willmore-sys}, the matrix $\mathbb{E}$ is
given by
\begin{equation}
\left\{
\begin{array}{l}
\displaystyle{\mathbb{E} = \frac{1}{Q_\varepsilon} \left(\mathbb{I}-\frac{\r \otimes
    \r }{Q_\varepsilon^2}\right),}
\\
\,
\\
Q_\varepsilon = \sqrt{\varepsilon+|\r|^2}.
\label{eq:eq}
\end{array}\right.
\end{equation}

Next, we obtain the weak formulation of the exact solution. Applying
the LDG method to system \eqref{eq:willmore-sys}, we get the following
numerical scheme for the unknowns $(\phi,\omega, H) \in V_h^3$ and  $(\s,\vn,\p,\q,\r) \in \Sigma_h^5$, such that, for all test functions $(\varphi,\xi, \vartheta) \in V_h^3$ and $(\thetan,\psin,\etan,\rhon,\zetan) \in \Sigma_h^5$, we have
\begin{equation}
\int_K \frac{1}{Q_\varepsilon} \,\frac{\partial \phi}{\partial t} \,\varphi \,dK \,=\,\int_K  (\s-\vn)\,
\cdot \nabla \varphi \,dK \,\,-\,\,\int_{\partial K} \left(\widehat{\s \cdot
  \nun}\,-\,\widehat{\vn \cdot \nun}\right)\,\varphi\, ds,
\label{eq:willmore-ldg}
\end{equation}
where the previous algebraic relations between $(\s,\vn)$ and
$(\p,\r,\omega)$ are now written as
\begin{equation*}
\left\{
\begin{array}{l}
\displaystyle\int_K \s \cdot \thetan \,dK \,=\, \int_K \mathbb{E} \,\p \cdot \thetan \,dK,
\\
\,
\\
\displaystyle\int_K \vn \cdot \psin \,dK \,\,=\,\, \int_K \frac{1}{2} \,\frac{\omega^2}{Q_\varepsilon^3} \,\r \cdot \psin \,dK.
\end{array}\right.
\end{equation*}
Moreover, $\p$ is solution to
$$
\int_K \p \cdot \etan \,dK \,\,=\,\,-\int_K \omega \,\nabla \cdot \etan
\,dK\,\,+\,\,\int_{\partial K}\widehat{\omega} \nun \cdot \etan \,ds,
$$
where $\omega$  is now given with respect to $H$
$$
\int_K \omega \xi \,dK \,\,=\,\, \int_K  Q_\varepsilon \,H \,\xi\, dK,
$$
and $H$ is solution to
$$
\int_K H \vartheta \,dK\,\,=\,\,-\int_K \q \,\cdot \,\nabla \vartheta
\,dK\,\,+\,\,\int_{\partial K} \widehat{\q \cdot \nun} \vartheta \,ds,
$$
where $(\q,\r)$ are given  with respect to the main unknown $\phi$,
that is,
\begin{equation*}
\left\{
\begin{array}{l}
\displaystyle \int_K \q \cdot \rhon \,dK \,\,=\,\,\int_K
\frac{\r}{Q_\varepsilon}\, \cdot \rhon \,dK,
\\
\,
\\
\displaystyle \int_K \r \cdot \zetan \,dK\,\, =\,\,-\int_K \phi \nabla \cdot \zetan
\,dK\,\,+\,\,\int_{\partial K} \widehat{\phi} \nun \cdot \zetan\, ds,
\end{array}\right.
\end{equation*}
where $\mathbb{E}$ and $Q_\varepsilon$ are computed by \eqref{eq:eq}. Here $\widehat{\s}$, $\widehat{\vn}$, $\widehat{\omega}$, $\widehat{\q}$ and $\widehat{\phi}$ are the numerical fluxes. To complete the definition of the LDG method, we need to define these numerical fluxes, which are functions defined on the edges and should be designed based on different guiding principles for different PDEs to ensure stability.

Similar to the development in \cite{Xu_will}, it turns out that we can take the simple choices such that
\begin{equation}\label{flux}
\widehat{\s}|_e=\s_L,\;\;\;\; \widehat{\vn}|_e=\vn_L, \;\;\;\; \widehat{\q}|_e=\q_R,\;\;\;\; \widehat{\omega}|_e=\omega_L,\;\;\;\; \widehat{\phi}|_e=\phi_R.
\end{equation}
We remark that the choice for the numerical fluxes \eqref{flux} is not unique. In fact, the crucial part is taking $\widehat{\s}$, $\widehat{\vn}$ and $\widehat{\phi}$ from opposite sides, and $\widehat{\q}$ and $\widehat{\omega}$ from opposite sides.

\begin{rem}
The solution to the LDG scheme  is mass conservative, by choosing the test function $\varphi=1$ in equation \eqref{eq:willmore-ldg}.
\end{rem}
\begin{rem}
Similar to the proof in \cite{Xu_will}, the solution to the LDG scheme  \eqref{eq:willmore-ldg} with numerical fluxes \eqref{flux} satisfies the following energy stability
\begin{equation*}
\frac{1}{2} \frac{d}{dt} \int_\Omega H^2 Q_\varepsilon d \x+\int_\Omega \frac{(\phi_t)^2}{Q_\varepsilon} d \x=0.
\end{equation*}
\end{rem}
\subsection{ A $p$-adaptive technique}
In the level set method, we aim to capture the interface behavior. As a result, we are mainly interested in the position near the zero level set, which motivates us to use adaptive technique. The basis functions of the LDG method we discuss in this paper can be discontinuous, so it has certain adaptive flexibility, such as, the order of accuracy can be locally determined in each cell.

\begin{figure}[!ht]\centering
\begin{tikzpicture}
\fill[gray!20](0,0) rectangle (1,1);
\fill[gray!20](1,0) rectangle (2,1);
\fill[gray!20](1,1) rectangle (2,2);
\fill[gray!20](1,-1) rectangle (2,0);
\fill[gray!20](0,-1) rectangle (1,0);
\fill[gray!20](1,-2) rectangle (2,-1);
\fill[gray!20](2,-2) rectangle (3,-1);
\fill[gray!20](1,-3) rectangle (2,-2);
\fill[gray!20](2,-3) rectangle (3,-2);
\fill[gray!20](-3,2) rectangle (-2,3);
\fill[gray!20](-2,2) rectangle (-1,3);
\fill[gray!20](-2,3) rectangle (-1,4);
\fill[gray!20](-2,1) rectangle (-1,2);
\fill[gray!20](-1,1) rectangle (0,2);
\fill[gray!20](0,1) rectangle (1,2);
\fill[gray!20](0,2) rectangle (1,3);
\fill[gray!20](-1,2) rectangle (0,3);
\fill[gray!20](0,-2) rectangle (1,-1);
\fill[gray!20](2,-1) rectangle (3,0);
\fill[gray!20](2,0) rectangle (3,1);
\fill[gray!20](-1,0) rectangle (0,1);
\fill[gray!20](-3,1) rectangle (-2,2);
\fill[gray!20](-3,3) rectangle (-2,4);
\fill[gray!20](-1,3) rectangle (0,4);
\fill[gray!20](-2,0) rectangle (-1,1);
\draw[step=1.0cm,black,thick] (-5,-5) grid (5,5);
\draw[thick] (2,-2.5) arc (0:75:5cm);
\end{tikzpicture}
\caption{Explanation of $p$-adaptive technique.}
\label{grid_adap}
\end{figure}
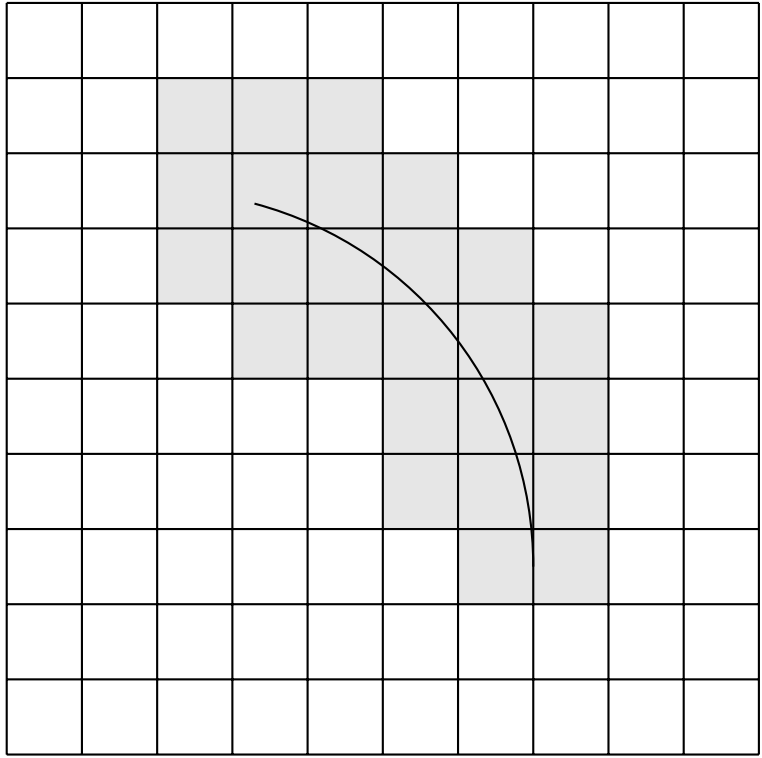

Our aim is to apply a  higher order polynomial approximation around
the zero level set to achieve high order accuracy, whereas in the
regions far away from the zero level set, we apply lower order polynomial approximations, which is similar as what is described by Yan and Osher in \cite{Yan}. We explain the $p$-adaptive technique in detail with Figure \ref{grid_adap}. In Figure \ref{grid_adap}, we apply $\mathcal{P}^2$ approximation near the zero level set, which are the domain filled with dark cells. While for the domain two or three cells away from the zero level set, it is not necessary to use high order approximation, thus  $\mathcal{P}^1$ approximation is enough.

As we know, in the implementation of the $p$-adaptive technique, it is essential to employ a flag to mark ``near the zero level set''. By the definition of the level set method, it is reasonable to use the following flag:
\begin{eqnarray}
\left\{
\begin{aligned}
 & \mathcal{P}^2 \;\; \text{approximation}, \;\; \text{if} \;\; \abs{\phi}\leq 1.5 h, \\
 &\mathcal{P}^1 \;\; \text{approximation}, \;\; \text{if} \;\; \abs{\phi} > 1.5 h,
  \end{aligned}
  \right .
\end{eqnarray}
where $\abs{\phi}$ is the absolute value of $\phi$. With the proposed $p$-adaptive technique, we can improve the efficiency of the local discontinuous Galerkin method for level set equation.

\section{Time discretization}
 \label{se:time}
 As we know, the Willmore flow \eqref{willmore} is a highly nonlinear fourth order PDE, explicit time integration methods lead to an extremely small time step restriction for stability, but not for accuracy. It would therefore be desirable to develop a semi-implicit time discretization technique to alleviate this problem.

In this section, we employ a high order semi-implicit Runge-Kutta method. Based on the work in \cite{Filbet}, the main idea of the semi-implicit method is to apply two different Runge-Kutta methods.
Consider the ODE system
\begin{equation} \label{ODE}
\left\{
\begin{aligned}
 &\frac{d u}{dt}\,=\,\mathcal{H}(t,u(t),u(t)), \;\; t \in[0,T]\\
 &u(0)=u_0,
  \end{aligned}
  \right .
  \end{equation}
where the function $\mathcal{H}:$ $\mathbb{R} \times \mathbb{R}^m \times \mathbb{R}^m \rightarrow \mathbb{R}^m$ is sufficiently differentiable and the dependence on the second argument of $\mathcal{H}$ is non-stiff, while the dependence on the third argument is stiff.

The semi-implicit Runge-Kutta method is based on the partitioned Runge-Kutta method, so we will first give a simple description for the partitioned Runge-Kutta methods.  We consider
autonomous differential equations in the partitioned form,
\begin{equation} \label{ODE_partion}
\left\{
\begin{aligned}
 &\frac{dy}{dt}(t)=\mathcal{F}(y(t),z(t)), \\
 &\frac{dz}{dt}(t)=\mathcal{G}(y(t),z(t)),
  \end{aligned}
  \right .
  \end{equation}
 where $y(t) \in \mathbb{R}^m$, $z(t)\in \mathbb{R}^n$, $\mathcal{F}: \mathbb{R}^m \times \mathbb{R}^n\rightarrow \mathbb{R}^m$, $\mathcal{G}:\mathbb{R}^m \times \mathbb{R}^n\rightarrow \mathbb{R}^n$ and $\mathcal{F}$, $\mathcal{G} \in \mathcal{C}^1(\mathbb{R}^m \times \mathbb{R}^n)$. $y(t_0)=y_0$, $z(t_0)=z_0$ are the initial conditions.  Then we can express the partitioned Runge-Kutta methods by applying two different Runge-Kutta methods as the following Butcher tableau:
 \[
 \begin{array}{c|c}
 \hat{c} & \hat{A}\\
 \hline
 &\hat{b}^{T}
 \end{array}\;\;\;\;\;\;\;\;\;\;\;\;\;\;\;
 \begin{array}{c|c}
 c & A\\
 \hline
 &b^{T}
 \end{array}
 \]

For practical reasons, in order to simplify the computations, we
consider that $\hat{A}$ is a strictly lower triangular matrix and
$A$ is a lower triangular matrix for the implicit part. In addition, the coefficients satisfy:
 \begin{equation}
\label{hyp:00}
 \hat{c}_i=\sum_{j=1}^{i-1}\hat{a}_{i,j}, \;\;\;\; \text{and} \;\;\;\; c_i=\sum_{j=1}^{i}a_{i,j}, \;\;\;\; \text{for} \;\;\;\; 1\leqslant i \leqslant s.
 \end{equation}
 By applying the partitioned Runge-Kutta time marching method, the
 solution of the autonomous system \eqref{ODE_partion} advanced from time $t^n$ to $t^{n+1}=t^n+\Delta t$ is given by
 \begin{equation}
\label{partition_kl}
\left\{
\begin{aligned}
 &k_i=\mathcal{F}\left(y^n\,+\,\Delta t \,\sum_{j=1}^{i-1}\hat{a}_{i,j}\,k_j,\,z^n\,+\,\Delta t \,\sum_{j=1}^i a_{i,j}\,l_j\right), \;\;\; 1\leqslant i \leqslant s,\\
 &l_i=\mathcal{G}\left(y^n\,+\,\Delta t \,\sum_{j=1}^{i-1}\hat{a}_{i,j}\,k_j,\,z^n\,+\,\Delta t \,\sum_{j=1}^i a_{i,j}\,l_j\right), \;\;\; 1\leqslant i \leqslant s,
  \end{aligned}
  \right .
  \end{equation}
  and we can calculate $y^{n+1}$ and $z^{n+1}$ as follows
  \begin{equation} \label{partition_uv}
\left\{
\begin{aligned}
 &y^{n+1}\,=\,y^n\,+\,\Delta t\,\sum_{i=1}^{s}\hat{b}_i\,k_i,\\
 &z^{n+1}\,=\,z^n\,+\,\Delta t\,\sum_{i=1}^{s}b_i\,l_i.
  \end{aligned}
  \right .
  \end{equation}

 To derive a semi-implicit Runge-Kutta
 scheme, we first rewrite the non autonomous differential equation
 \eqref{ODE} as an autonomous system where we double the
 number of variable, that is,
\begin{equation} \label{ODE_rewrite}
\left\{
\begin{aligned}
 &\frac{d}{dt}\left(\begin{array}{l} t \\ u(t) \end{array} \right) \,= \,\left(\begin{array}{l} 1 \\ \mathcal{H}(t,u(t),v(t)) \end{array} \right),\\
 &\frac{dv(t)}{dt} \;=\,\mathcal{H}(t,u(t),v(t)),
\\
&u(0)\,=\,v(0)\,=\,u_0
  \end{aligned}
  \right .
  \end{equation}
By a uniqueness argument the
solution to (\ref{ODE_rewrite}) corresponds to the one of (\ref{ODE}),
but now this system is written as an autonomous partitioned system (\ref{ODE_partion}), with  $y(t)=(t,u(t))$, $\mathcal{F}=(1,\mathcal{H})$ and
 $z(t)=u(t)$, $\mathcal{G}=\mathcal{H}$ and $y(t_0)=(t_0,u_0)$,
 $z(t_0)=u_0$ are the initial conditions.

Applying the partitioned Runge-Kutta scheme
  \eqref{partition_kl}-\eqref{partition_uv} to system
  \eqref{ODE_rewrite}, we can get a high order semi-implicit
  Runge-Kutta method for \eqref{ODE} : the first component of
  the first equation (\ref{ODE_rewrite}) only gives
$$
\hat c_i = \sum_{j=1}^i \hat a_{i,j},
$$
whereas the second component of the first equation and the second
equation of (\ref{ODE_rewrite}) are identical, which gives the
following semi-implicit scheme
$$
k_i = \mathcal{H}\left(t^n+ \hat c_i \Delta t, u^n\,+\,\Delta t \,\sum_{j=1}^{i-1}\hat{a}_{i,j}\,k_j,u^n\,+\,\Delta t \,\sum_{j=1}^i a_{i,j}\,k_j\right), \;\;\; 1\leqslant i \leqslant s.
$$
It is worth to mention here that $k_i$ is defined implicitly since
$a_{i,i}\neq 0$. Therefore, starting from $u^n$, we give the algorithm to calculate $u^{n+1}$ in the following.
  \begin{enumerate}
  \item Set  for $i=1,\ldots,s$,
  \begin{align}
\label{uv}
  &U_i=u^n+\Delta t \sum_{j=1}^{i-1}\hat{a}_{i,j}k_j, \nonumber \\
  &V_i=u^n+\Delta t \sum_{j=1}^{i}a_{i,j}k_j.
  \end{align}
  \item For $i=1,\ldots,s$, compute
  \begin{equation}\label{k}
  k_i=\mathcal{H}(t^n+\hat{c}_i \Delta t, U_i,V_i).
  \end{equation}
  \item Update the numerical solution $u^{n+1}$ as
  \begin{equation}
\label{algorithm1}
u^{n+1}=u^n+\Delta t \sum_{i=1}^{s}b_i \,k_i.
\end{equation}
  \end{enumerate}

After the LDG spatial discretization for the Willmore flow \eqref{willmore}, we can apply the semi-implicit scheme \eqref{uv}-\eqref{algorithm1} by writing the system of ODE in the form of \eqref{ODE} with the component $u(t)$ treated explicitly, the component $v(t)$ treated implicitly and
$$
\left\{
\begin{array}{l}
\mathcal{H}(t,u,v) \,\,=\,\,-Q_\varepsilon(u)\,\nabla \cdot
\left(\mathbb{E}(u) \nabla
  \omega(u,v)\,-\,\frac{1}{2}\frac{\omega(u)^2}{Q_\varepsilon(u)^3}\,
  \nabla v\right),
\\
\omega(u,v) \,\,=\,\,Q_\varepsilon(u) \,\nabla \cdot \left( \frac{\nabla v}{Q_\varepsilon(u)}\right),
\end{array}\right.
$$
where
$$
\left\{
\begin{array}{l}
Q_\varepsilon(u) \,\,=\,\,\sqrt{\varepsilon\,+\,|\nabla u|^2},\quad
\omega(u)\,=\,Q_\varepsilon(u) \,\nabla \cdot \left( \frac{\nabla
    u}{Q_\varepsilon(u)}\right),
\\
\mathbb{E}(u)\,\,=\,\,\frac{1}{Q_\varepsilon(u)}\,\left(\mathbf{I}\,-\,\frac{\nabla u}{Q_\varepsilon(u)}\otimes \frac{\nabla u}{Q_\varepsilon(u)} \right).
\end{array}\right.
$$

By algorithm \eqref{uv}-\eqref{algorithm1}, the second variable of
equation \eqref{ODE} is treated explicitly and the third one
is treated implicitly. Obviously, it is necessary to solve  system of
linear  algebraic equations \eqref{k} at each time
step. The overall performance highly depends on the performance of the
solver. Traditional iterative solution methods such as Gauss-Seidel
method suffers from slow convergence rates, especially for large
systems. To enhance the efficiency of the high order semi-implicit
time marching method, we will apply the multigrid solver  to solve
algebraic equations \eqref{k} in this paper. And for a detailed
description of the multigrid solver, we refer the readers to
Trottenberg {\it et al.} \cite{Tro} and Brandt \cite{achi:book}.
\section{Numerical results}
\label{se:numerical}
In this section, we present some numerical results for the Willmore flow. In addition, we use $p$-adaptive LDG spatial discretization ($\mathcal{P}^2$ approximation near the zero level set, $\mathcal{P}^1$ approximation away from the zero level set) coupled with a high order semi-implicit Runge-Kutta method, which is high order accurate in both time and space, allowing larger time step. Each time step we solve the linear algebraic equations by multigrid solver. We first present various representative examples of shape relaxation, for example, the motions of an ellipse, a square, an asteroid and a singular curve. Next, we will show the topological changes, such as merging and breaking.

\subsection{ Accuracy test}

In this example, we consider the accuracy test for one-dimensional
Willmore flow. We consider the exact solution
\begin{equation}\label{eq:accuracy}
\phi(x,t)=0.05\sin(x)\cos(t)
\end{equation}
for equation \eqref{willmore} with a source term $f(x,t)$, which is
a given function so that \eqref{eq:accuracy} is the exact solution.

{To show the effect of different regularization parameters $\varepsilon$ on the corresponding errors and orders of accuracy,} we report the $L^2$ and $L^\infty$ errors and the numerical orders of accuracy at time $T=0.5$ with uniform meshes in Table \ref{table:accuracy}. We can see that the method with $\mathcal{P}^k$ elements gives $(k+1)$-th order of accuracy with different scales of the regularization parameters $\varepsilon$.
\begin{table}[!ht]
 \centering
\bigskip
\begin{footnotesize}
\begin{tabular}{|c||c|c|c|c||c|c|c|c||c|c|c|c|}
  \hline
   &\multicolumn{4}{|c||} {$\varepsilon=1$}&\multicolumn{4}{|c||} {$\varepsilon=5 h$} & \multicolumn{4}{|c|}{$\varepsilon=5 h^2$}\\\hline
 Mesh & $L^2 $ error & order &$ L^\infty$ error & order&$L^2 $ error & order &$ L^\infty$ error & order& $L^2 $ error & order &$ L^\infty$ error & order \\\hline
\multicolumn{13}{|c|}{$\mathcal{P}^1$}\\\hline
 16 &  4.16E-03  &   --  &  3.52E-03  &  --   &    4.16E-03  &  --  &3.52E-03&  --&    4.16E-03  &  --  &3.52E-03&  --\\\hline
 32 &  1.03E-03  &  2.00 &  9.07E-04  &  1.96 &    1.03E-03  &  2.00&9.09E-04&  1.95&    1.03E-03  & 2.00  &9.06E-04&  1.96\\\hline
 64 &  2.58E-04  &  2.00 &  2.27E-04  &  1.99 &   2.58E-04  &  2.00&2.27E-04&  2.00&    2.59E-04  &  1.99 &2.19E-04&  2.04 \\\hline
 \multicolumn{13}{|c|}{$\mathcal{P}^2$}\\\hline
 16 &  2.60E-04  &   --  &  2.07E-04  &  --   &    2.60E-04  &   --  &  2.08E-04  &  --&   2.60E-04  &   --  &  2.08E-04  &  --\\\hline
 32 &  3.25E-05  &  3.00 &  2.69E-05  &  2.95 &    3.25E-05  &  3.00 &  2.69E-05  &  2.95&    3.25E-05  &  3.00 &  2.69E-05  &  2.95\\\hline
 64 &  4.06E-06  &  3.00 &  3.40E-06  &  2.99 &    4.06E-06  &  3.00 &  3.40E-06  &  2.99 &    4.06E-06  &  3.00 &  3.39E-06  &  2.99\\\hline
\end{tabular}
\end{footnotesize}
\caption{\label{table:accuracy} Accuracy test. $L^2$ and $L^{\infty}$ error
  norms on the solution   at time $T=0.5$ with different regularization parameters $\varepsilon$.}
\end{table}

\subsection{Towards a circle}
We first perform a series of numerical simulations, where the initial
level sets describe two convex curves, a non-convex curve and a singular curve. Hence it is expected
that the level set approaches a circle for large time, { which compares very well with numerical calculations performed by Bene\v{s} {\it et al.} \cite{Bene}.}

\subsubsection*{The motion of an ellipse}

\begin{figure}[!ht]
\centering
\subfigure[$t=0$] {
\includegraphics[width=0.4\columnwidth]{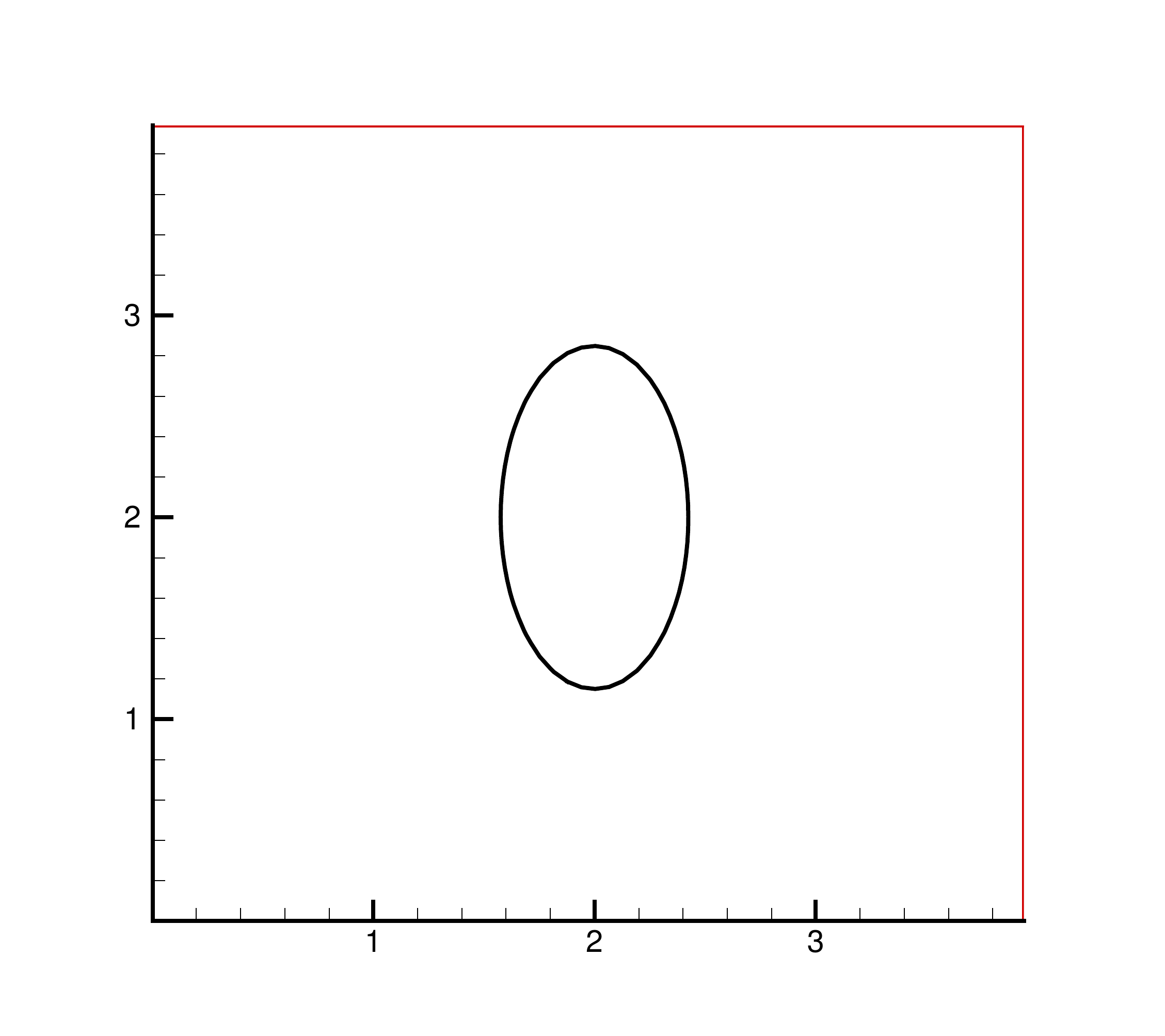}
}
\subfigure[$t=0.01$] {
\includegraphics[width=0.4\columnwidth]{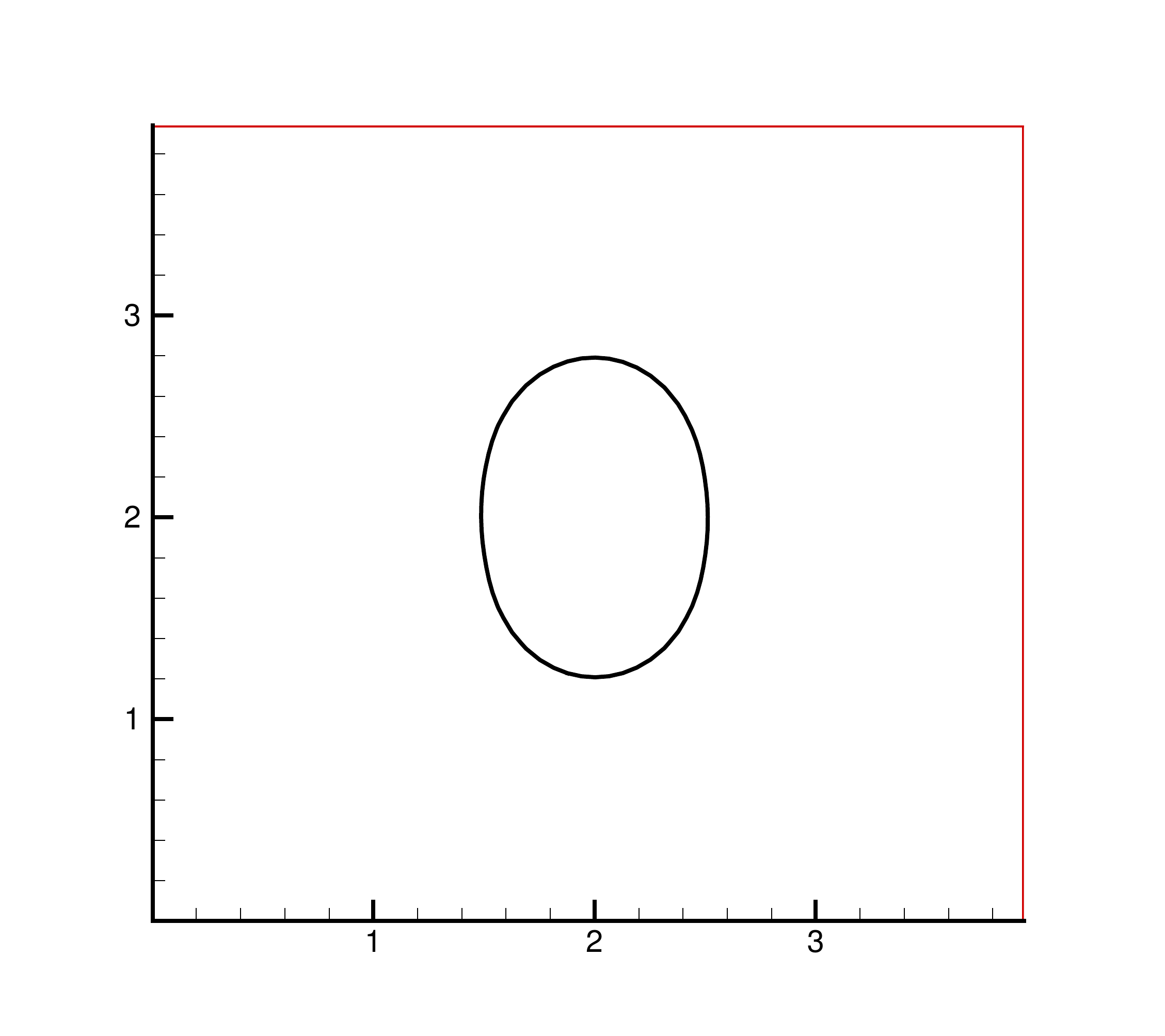}
}
\subfigure[$t=0.05$] {
\includegraphics[width=0.4\columnwidth]{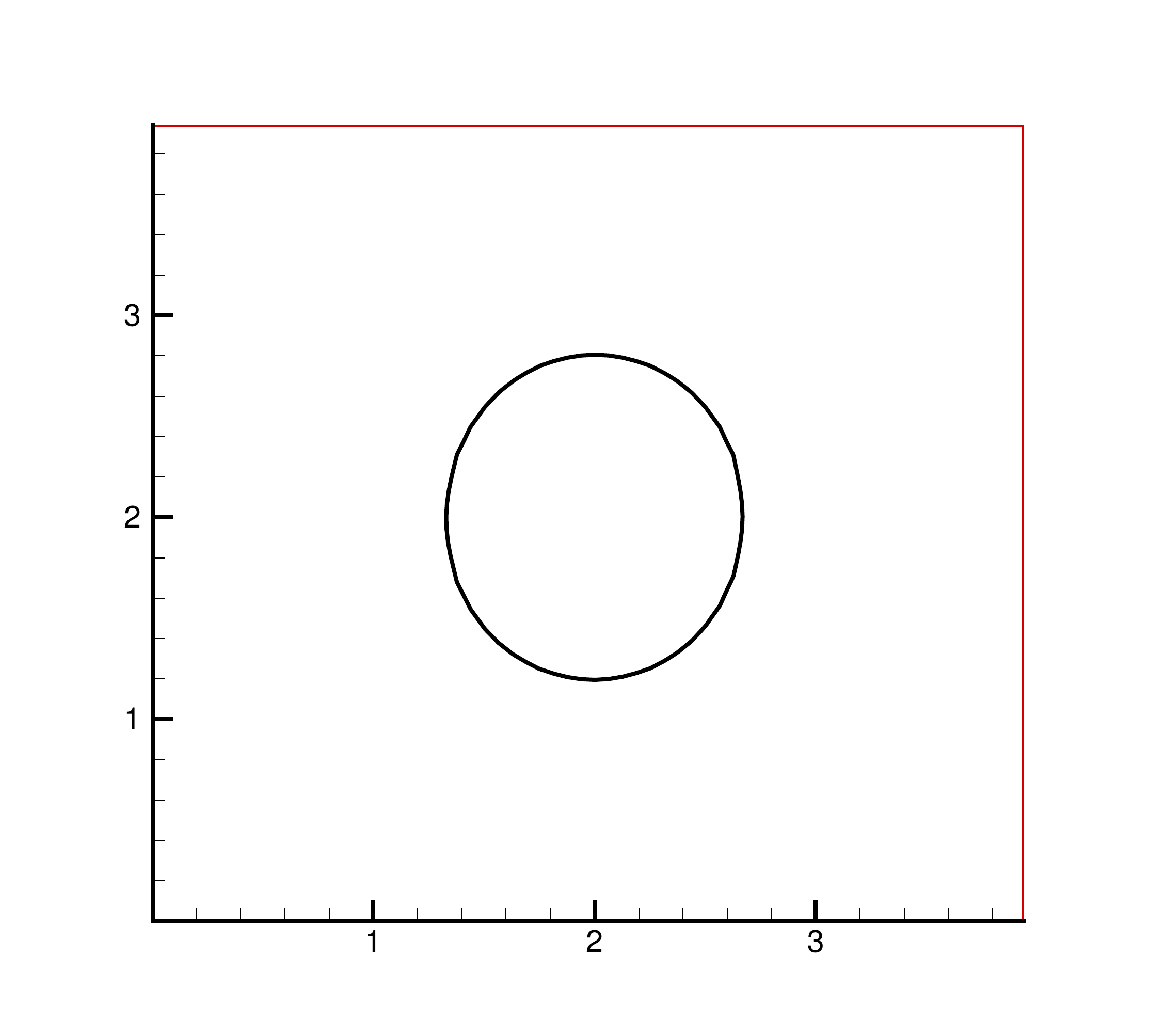}
}
\subfigure[$t=0.1$] {
\includegraphics[width=0.4\columnwidth]{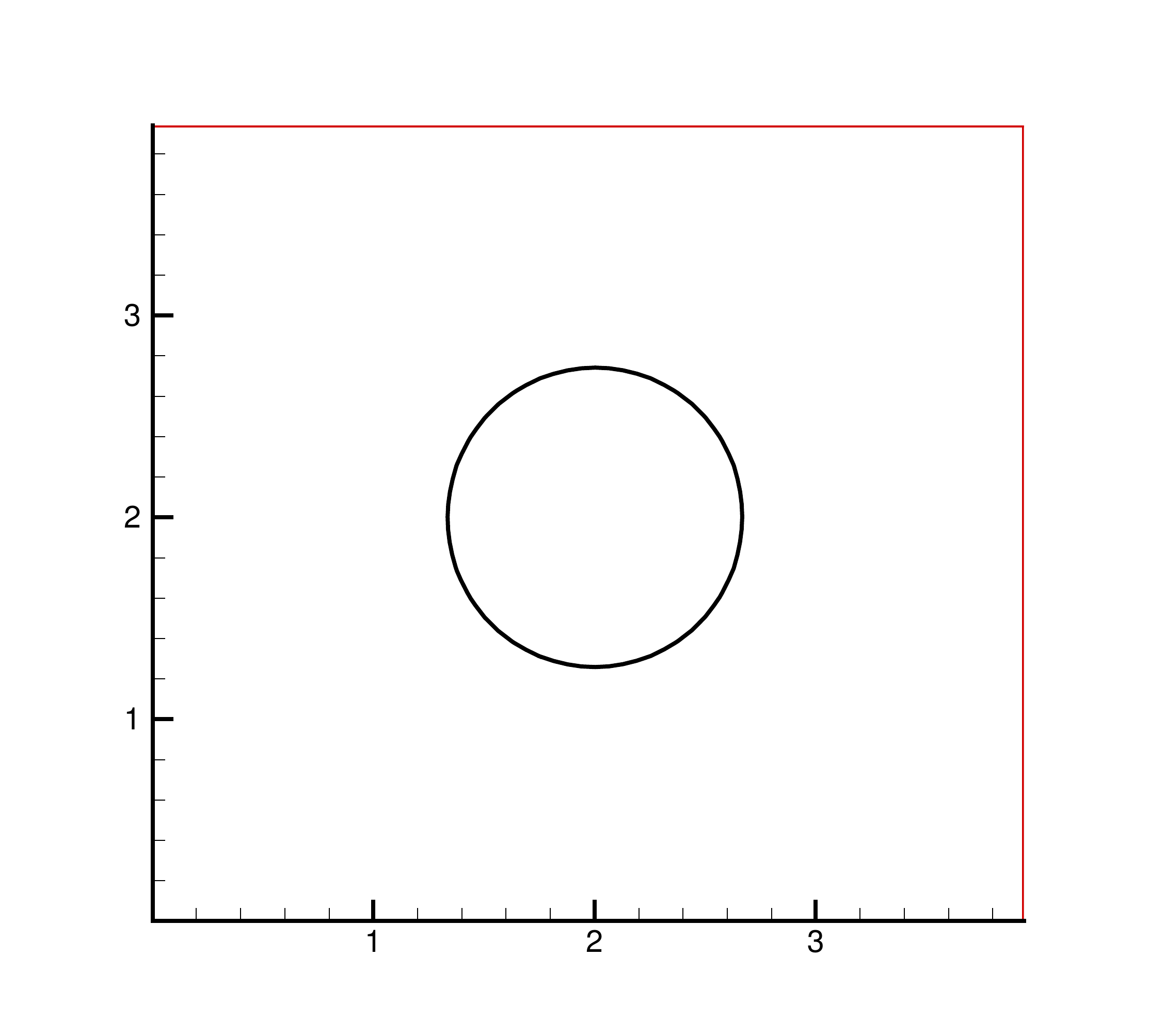}
}
\caption{Towards a circle. {\it Evolution of Willmore of the
    regularized  flow with $\varepsilon=h$ when the initial level set
    corresponds to an ellipse.}}
\label{ellipse}
\end{figure}

We consider the motion of an
ellipse as shown in Figure \ref{ellipse}. For large time,  the level set  asymptotically approaches a circle. Here, we approximate the Willmore flow with $\varepsilon=h$ on a $64 \times 64$ grid for the computational domain $\Omega=[0,4] \times [0,4]$. Here, we adopt a $p$-adaptive LDG method for spatial discretization and a second order semi-implicit Runge-Kutta method with a time step of $\Delta t=0.001 h$, which is very large compare to that of an explicit time marching method ($\Delta t=C h^4$).

\subsubsection*{The motion of a square}
\begin{figure}[!ht]
\centering
\subfigure[$t=0$] {
\includegraphics[width=0.4\columnwidth]{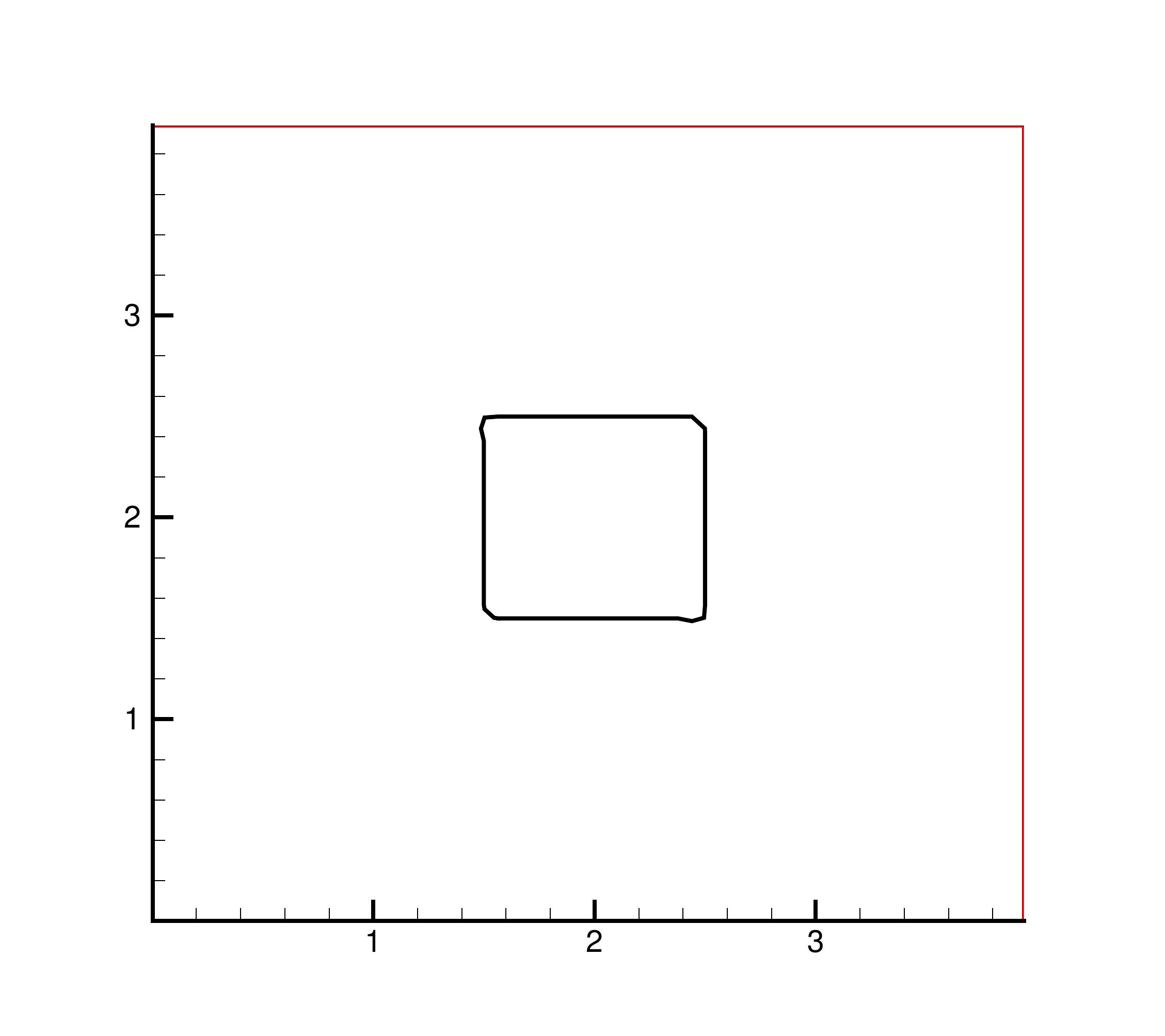}
}
\subfigure[$t=0.0005$] {
\includegraphics[width=0.4\columnwidth]{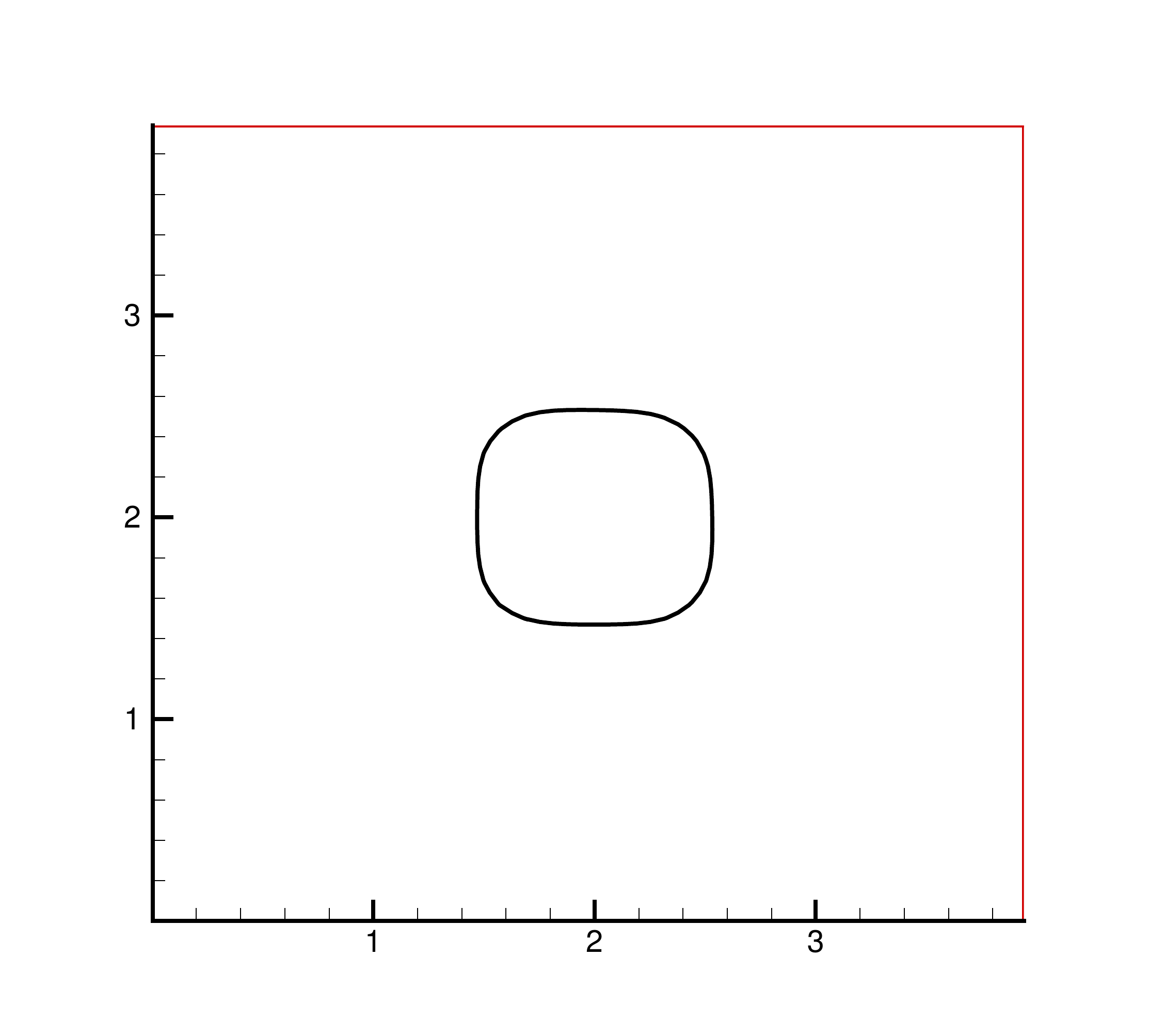}
}
\subfigure[$t=0.001$] {
\includegraphics[width=0.4\columnwidth]{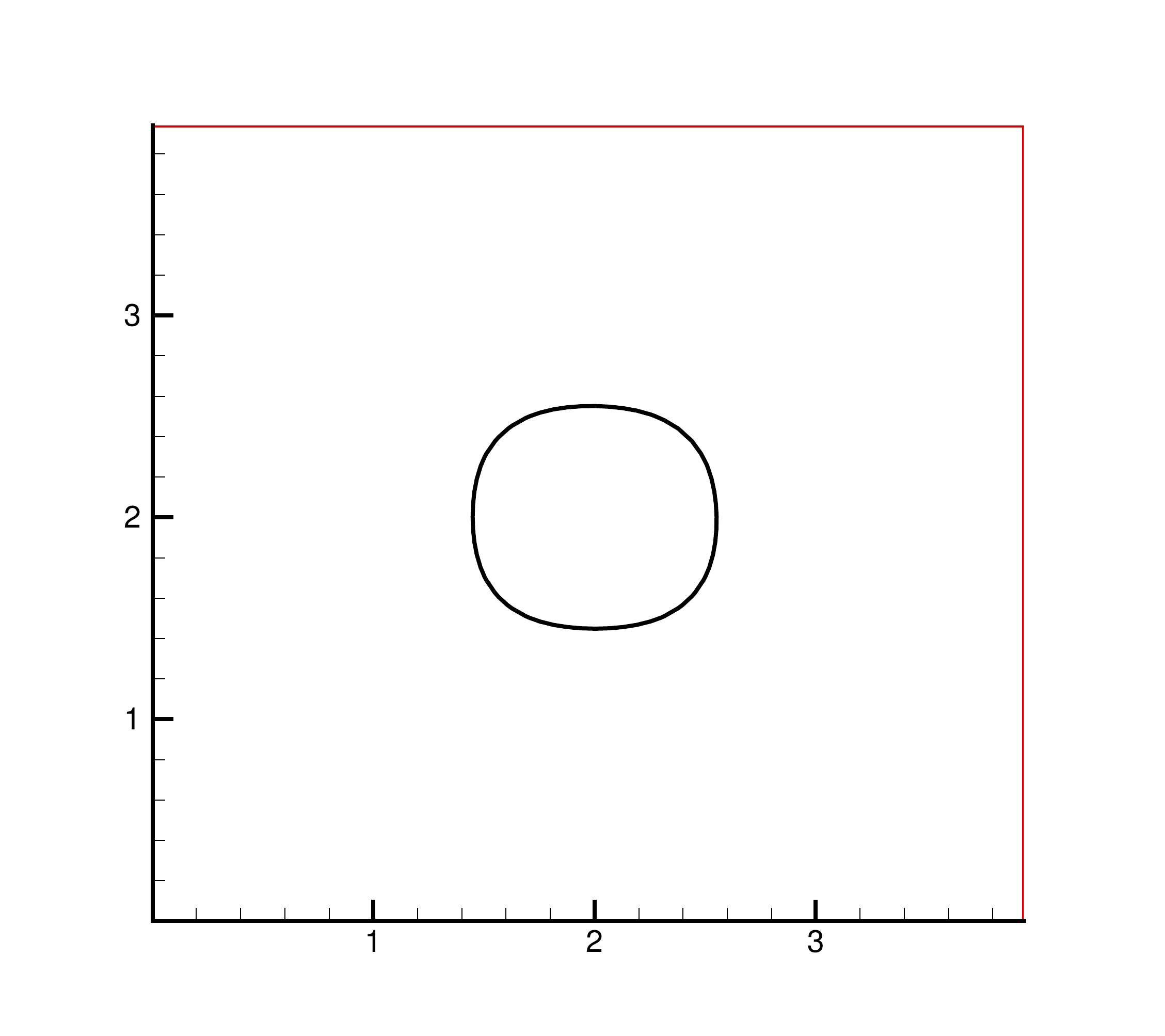}
}
\subfigure[$t=0.002$] {
\includegraphics[width=0.4\columnwidth]{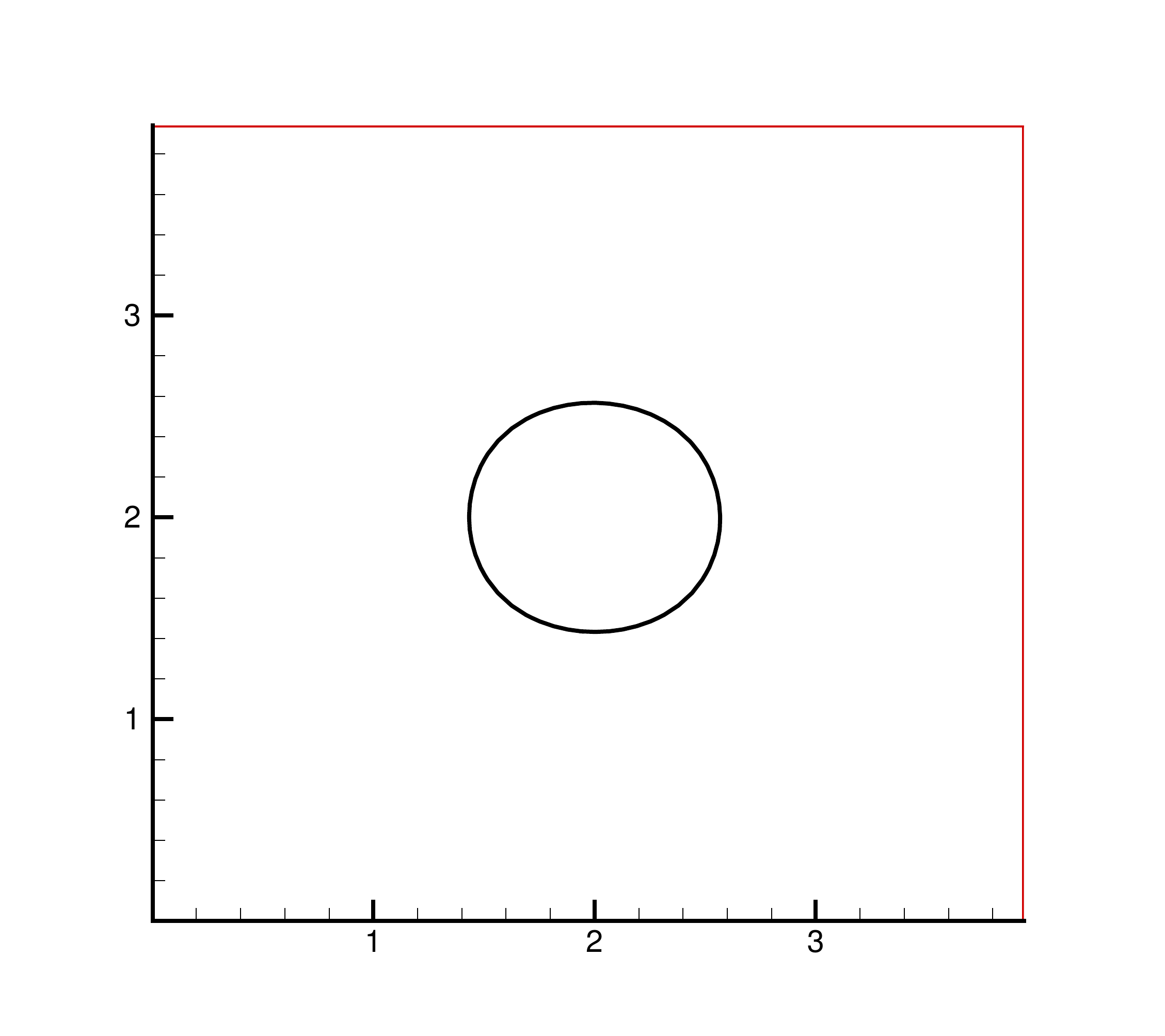}
}
\caption{Towards a circle.  {\it Evolution of Willmore of the
    regularized  flow with $\varepsilon=h$ when the initial level set
    corresponds to a square.}}
\label{square}
\end{figure}

We now consider an initial data, where the initial level set is a square,
that is, it has sharp corners.  We report the numerical results in
Figure \ref{square},  which shows  the evolution of a square into a circle. We approximate the Willmore flow with $\varepsilon=h$ on a $64 \times 64$ grid for the computational domain $\Omega=[0,4] \times [0,4]$. Here, we also adopt a second order semi-implicit Runge-Kutta method with a time step of $\Delta t=0.001 h$.

\subsubsection*{The motion of an asteroid}

\begin{figure}[!ht]
\centering
\subfigure[$t=0$] {
\includegraphics[width=0.4\columnwidth]{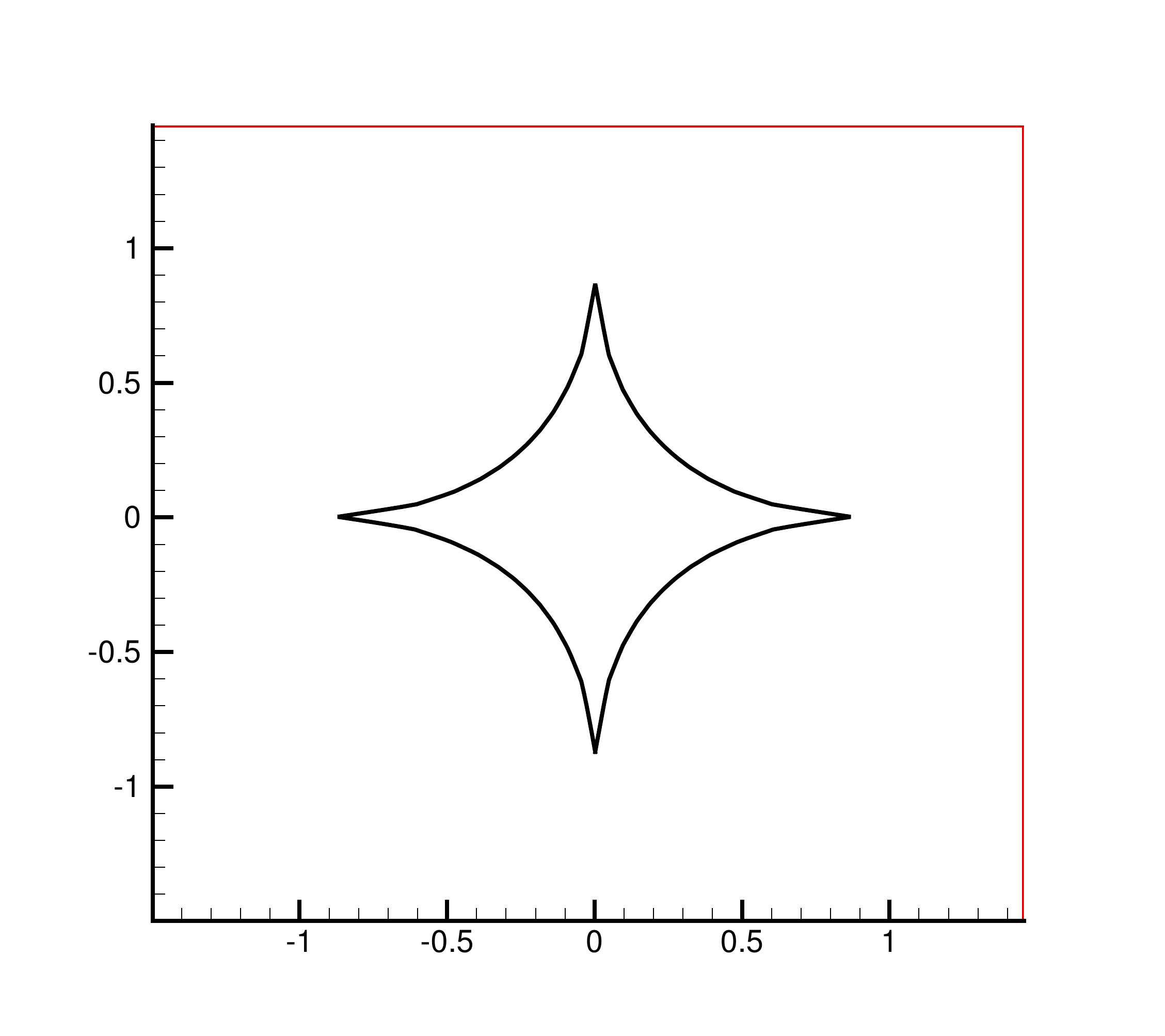}
}
\subfigure[$t=0.0001$] {
\includegraphics[width=0.4\columnwidth]{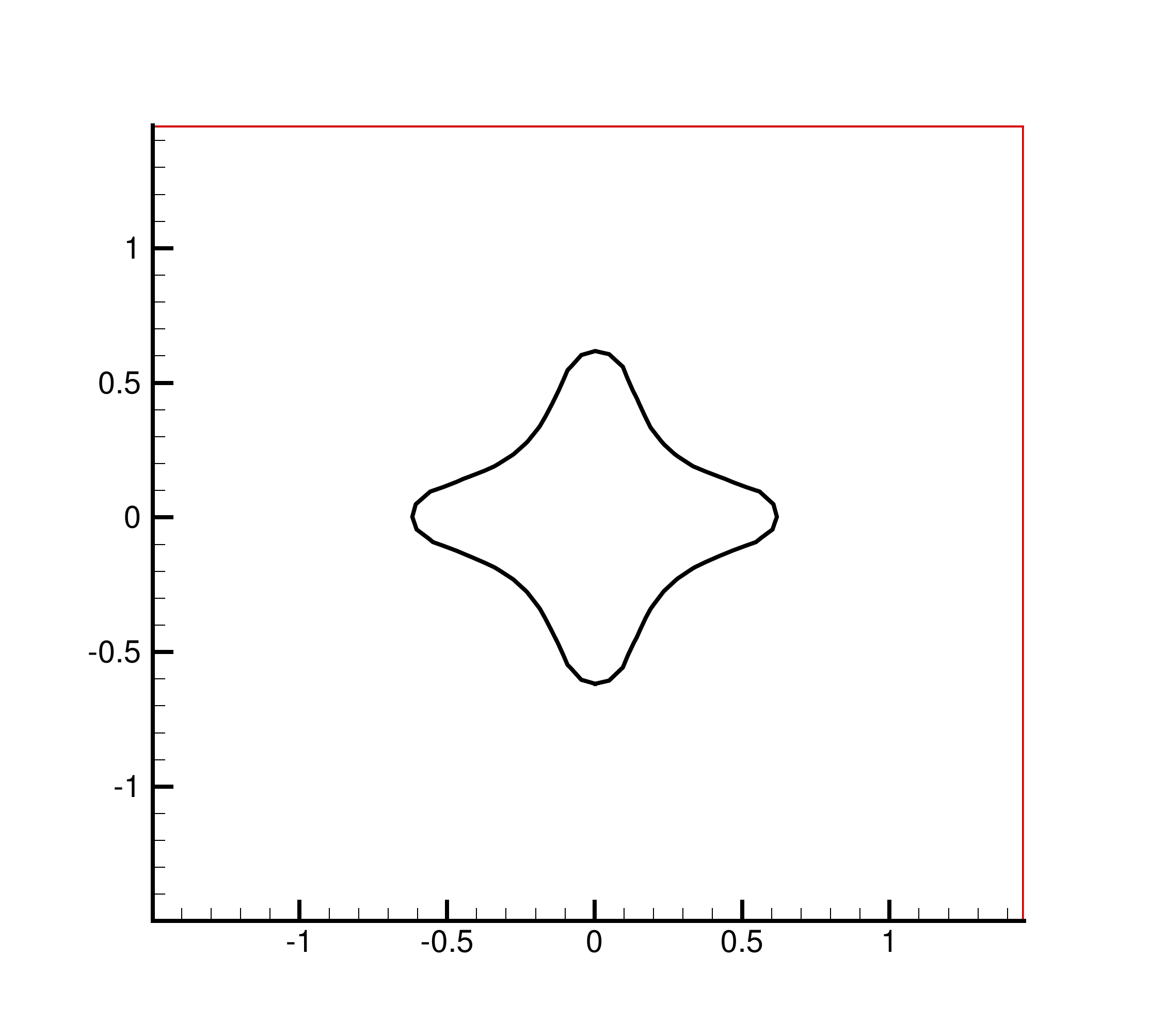}
}
\subfigure[$t=0.0005$] {
\includegraphics[width=0.4\columnwidth]{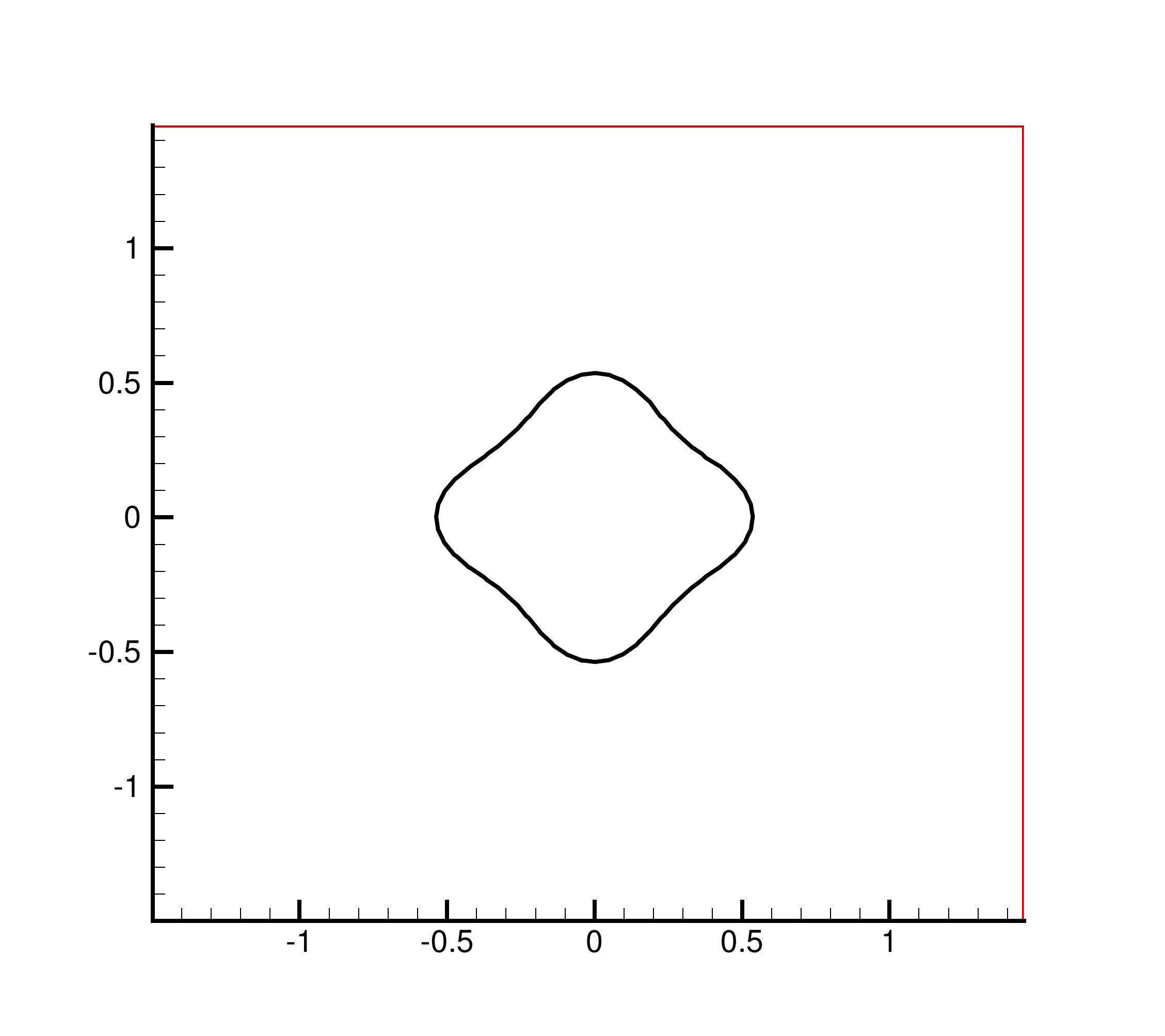}
}
\subfigure[$t=0.005$] {
\includegraphics[width=0.4\columnwidth]{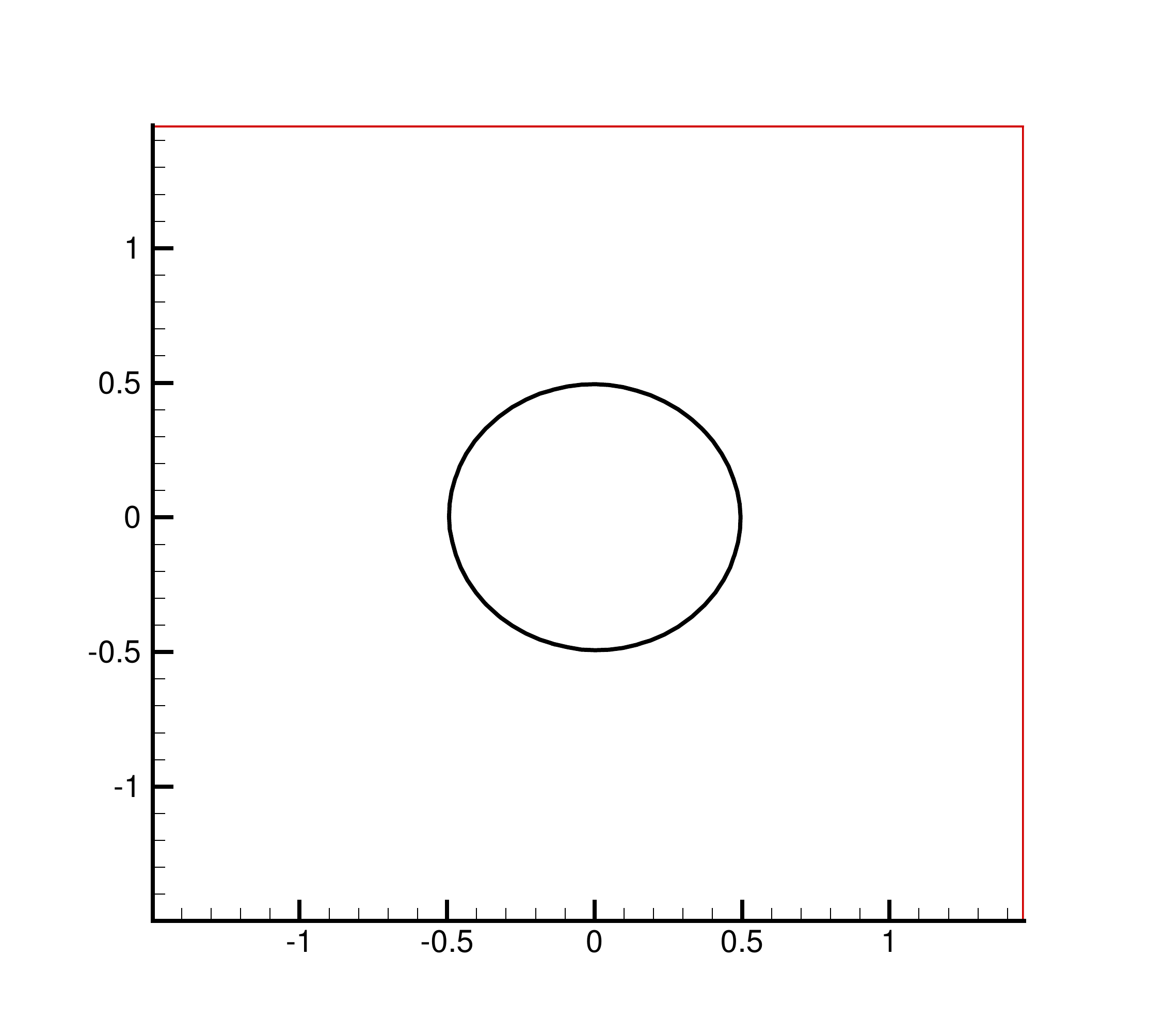}
}
\caption{Towards a circle.  {\it Evolution of Willmore of the
    regularized  flow with $\varepsilon=h$ when the initial level set
    corresponds to an asteroid.}}
\label{star}
\end{figure}

In Figure \ref{star}, we  present the motion of another non-convex curve with very sharp corners. The numerical methods and parameters  are chosen to be the same with those in the previous example. Also the initial asteroid is evolved to be a circle.

\subsubsection*{The motion of a singular curve}
\begin{figure}[!ht]
\centering
\subfigure[$\varepsilon=10h$, $t=0$] {
\includegraphics[width=0.32\columnwidth]{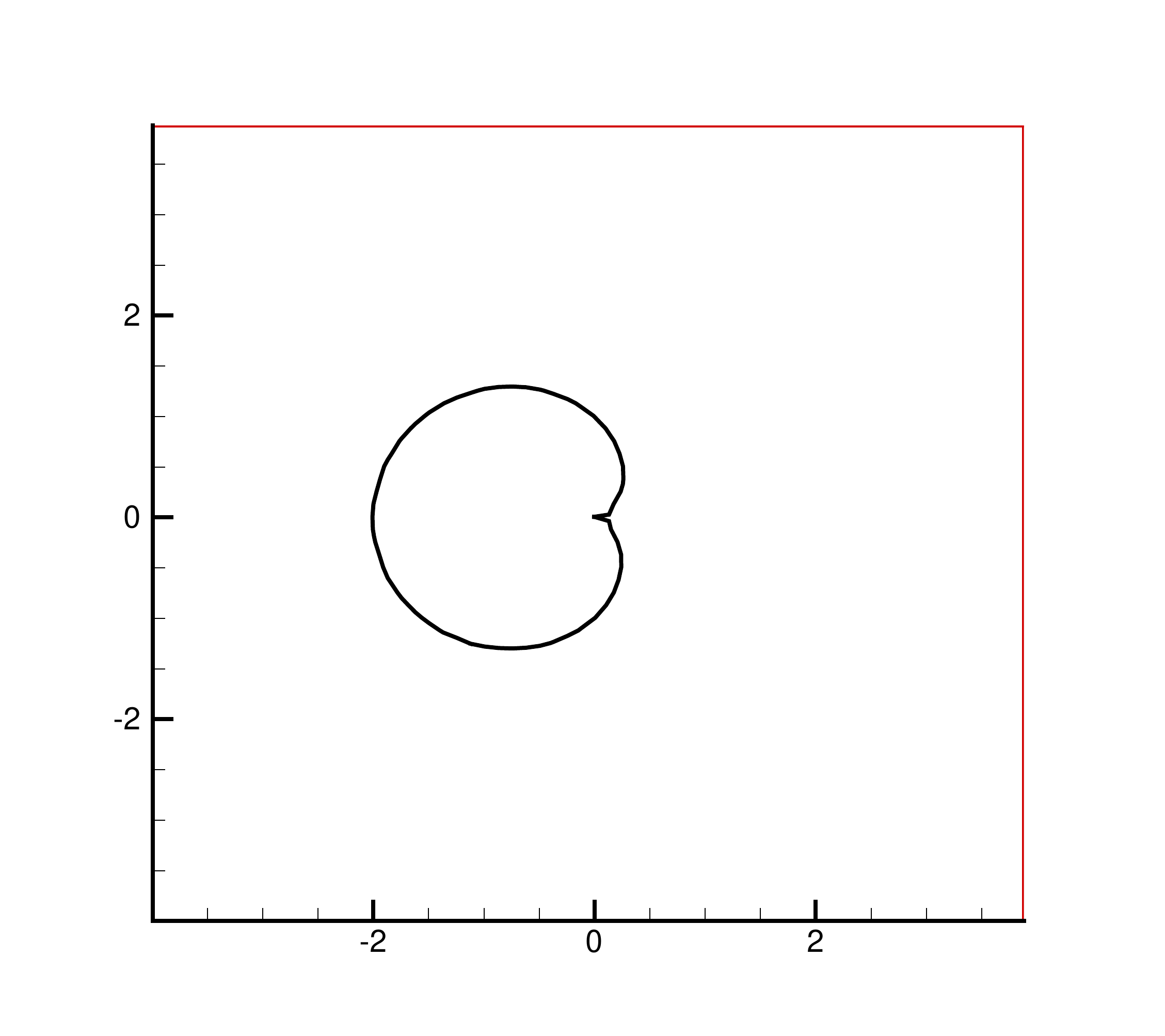}
}
\subfigure[$\varepsilon=10h^2$, $t=0$] {
\includegraphics[width=0.32\columnwidth]{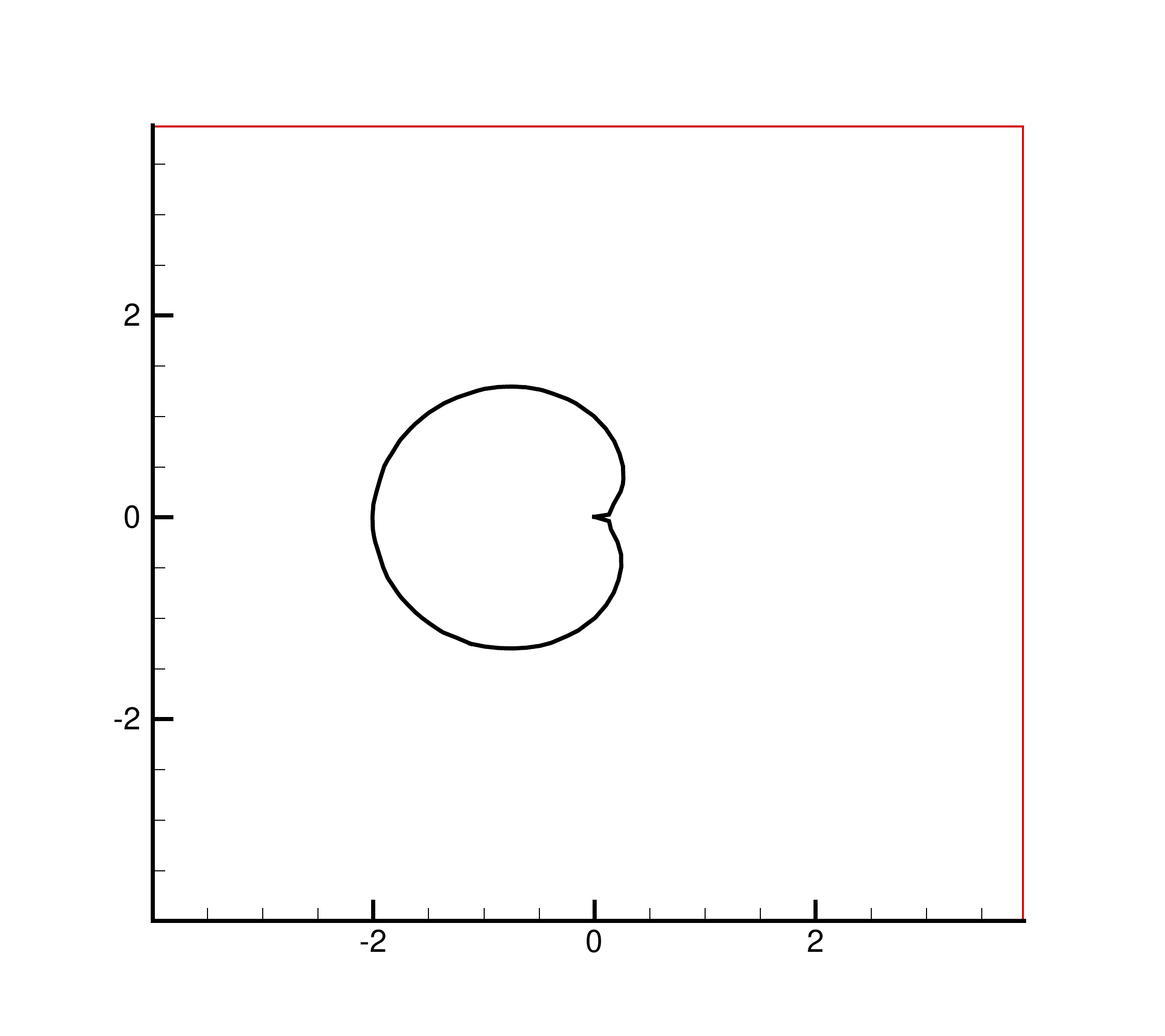}
}
\subfigure[$\varepsilon=10h$, $t=0.0005$] {
\includegraphics[width=0.32\columnwidth]{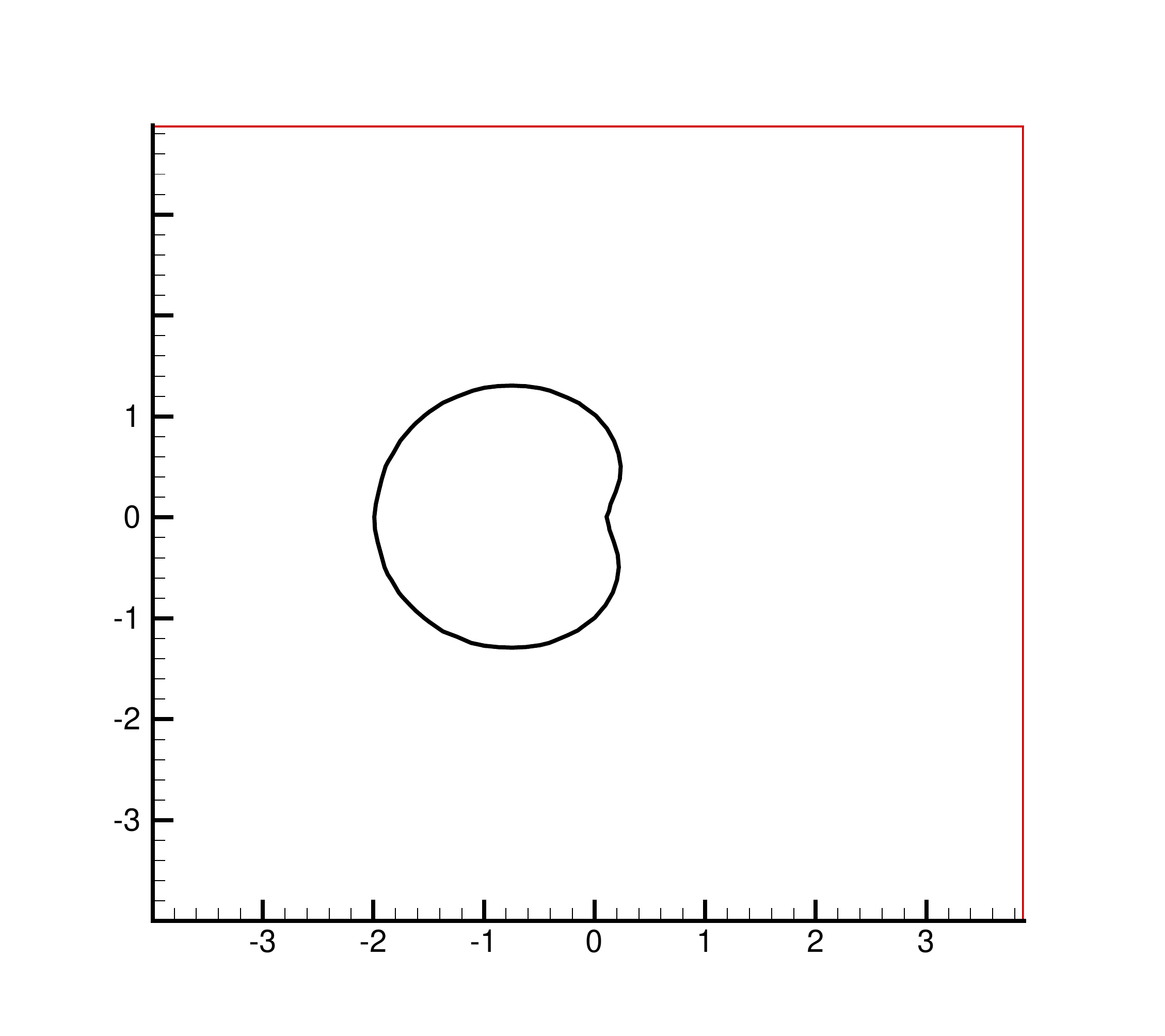}
}
\subfigure[$\varepsilon=10h^2$, $t=0.0005$] {
\includegraphics[width=0.32\columnwidth]{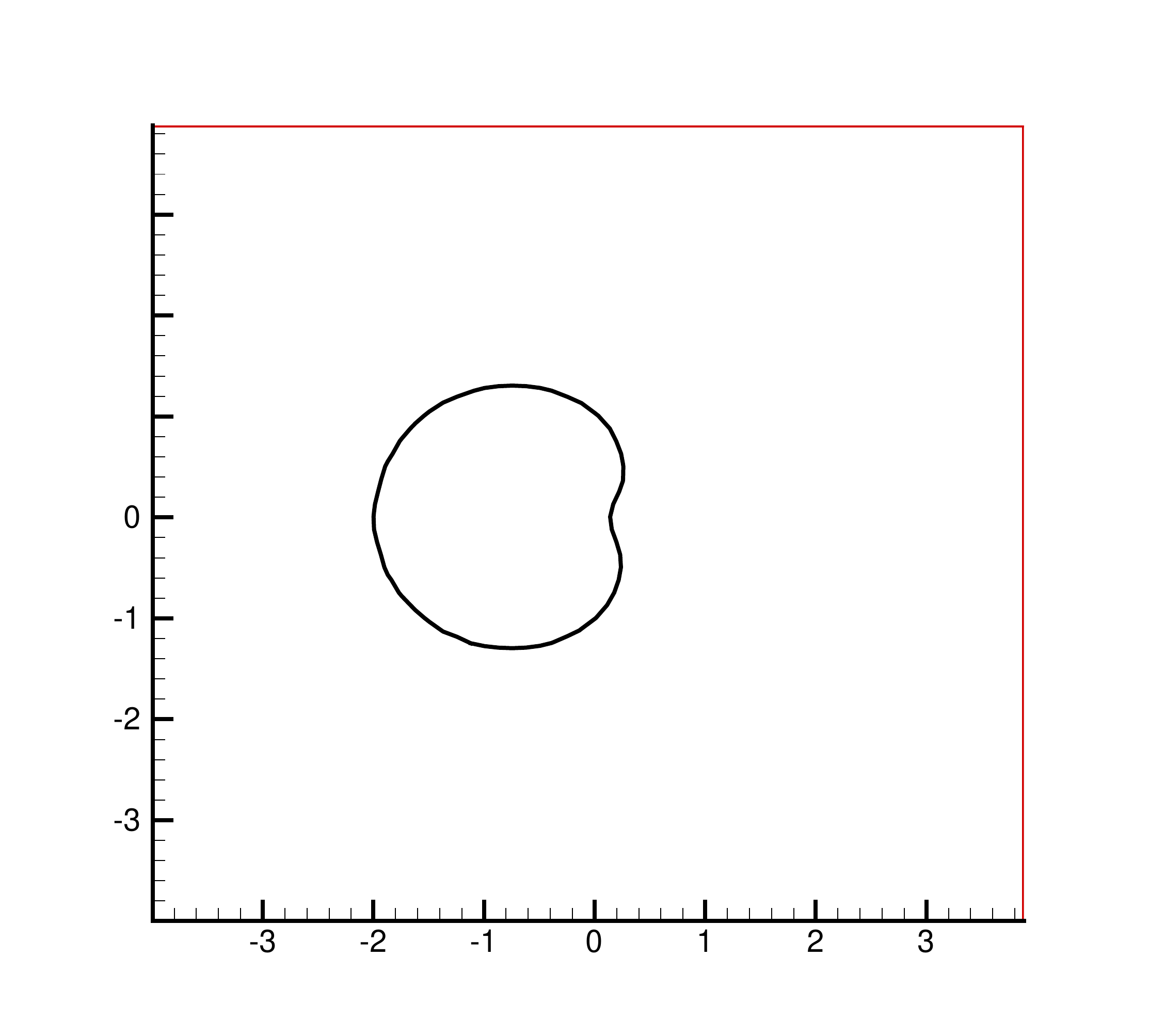}
}
\subfigure[$\varepsilon=10h$, $t=0.01$] {
\includegraphics[width=0.32\columnwidth]{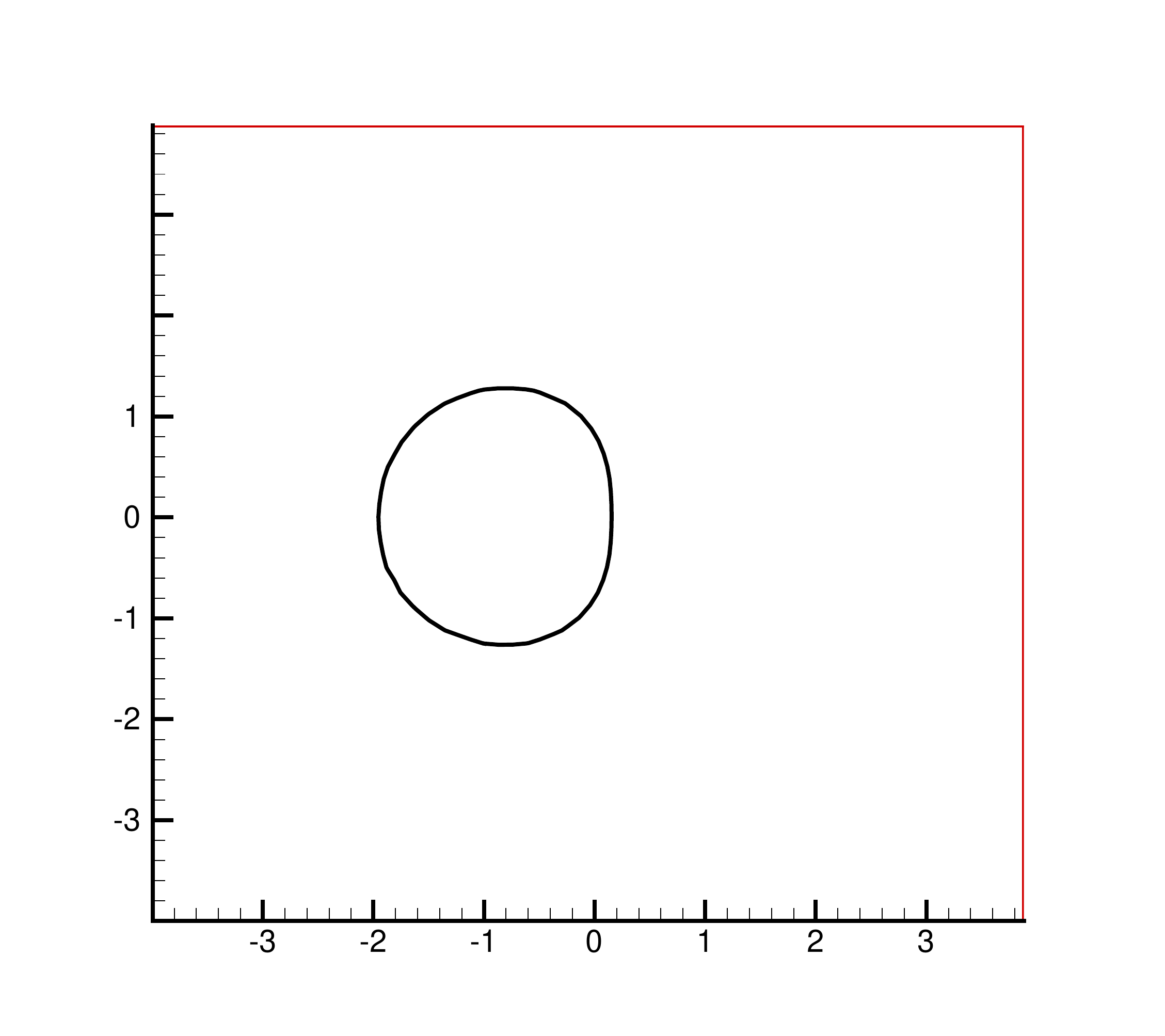}
}
\subfigure[$\varepsilon=10h^2$, $t=0.01$] {
\includegraphics[width=0.32\columnwidth]{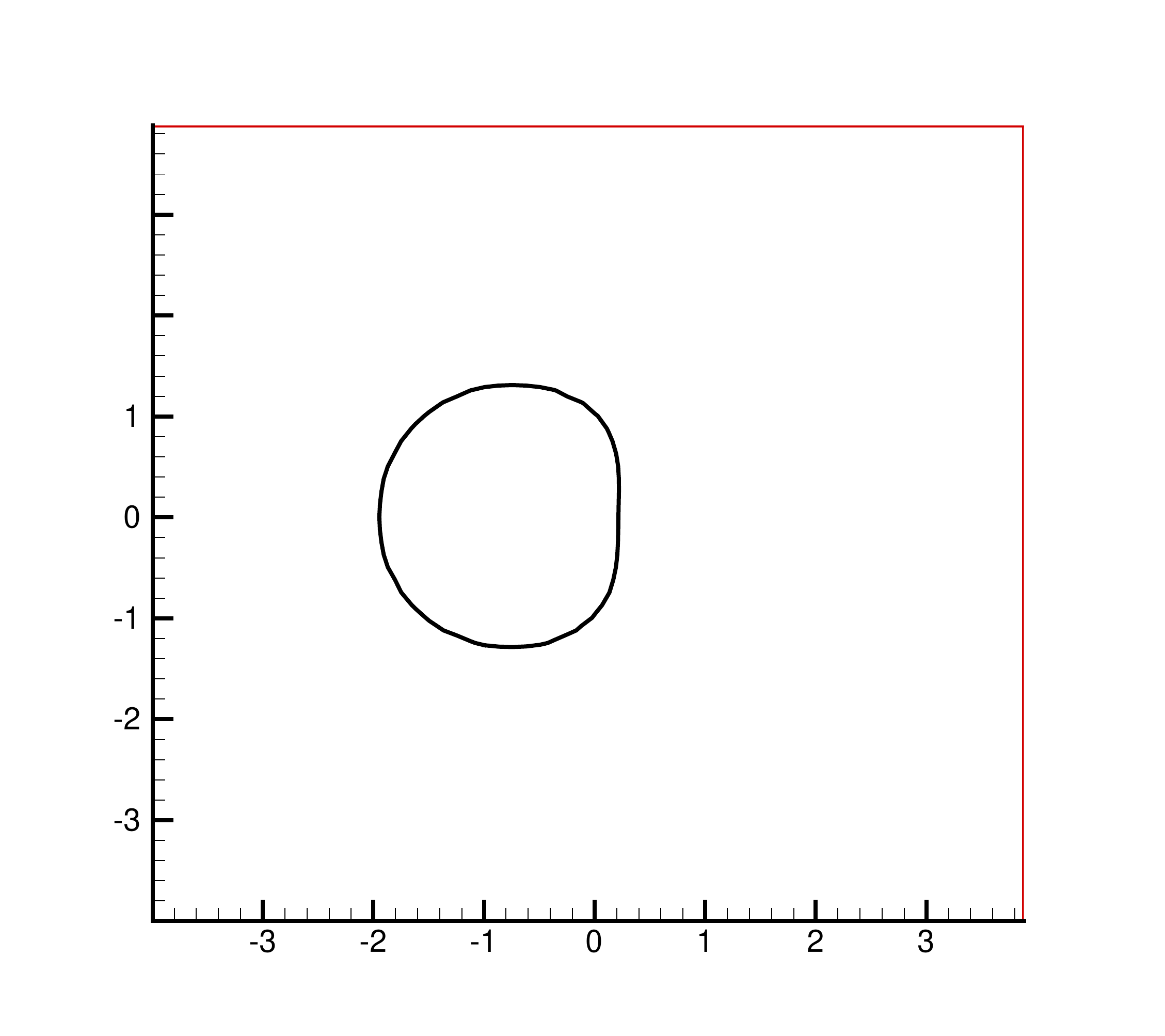}
}
\subfigure[$\varepsilon=10h$, $t=0.1$] {
\includegraphics[width=0.32\columnwidth]{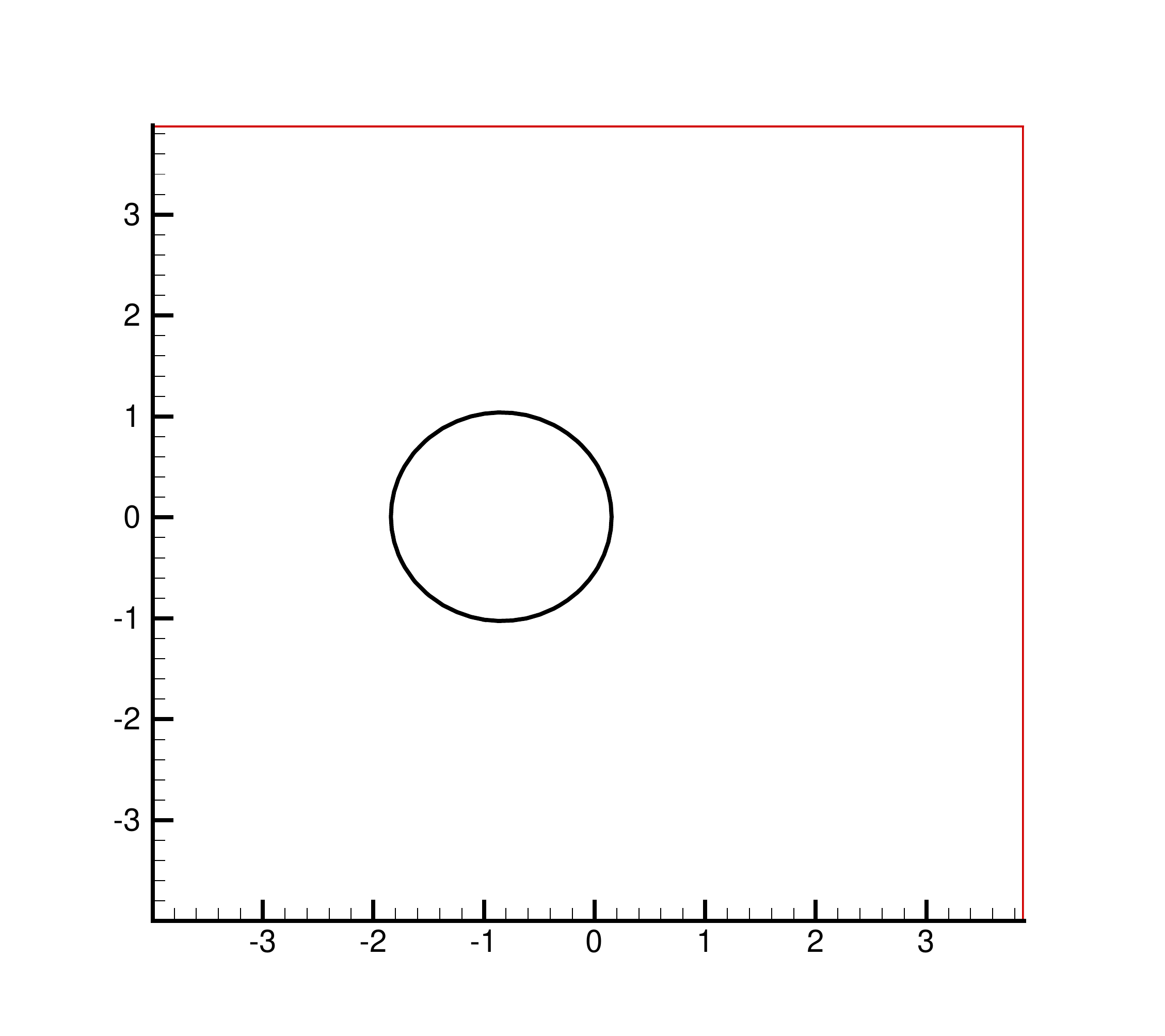}
}
\subfigure[$\varepsilon=10h^2$, $t=0.1$] {
\includegraphics[width=0.32\columnwidth]{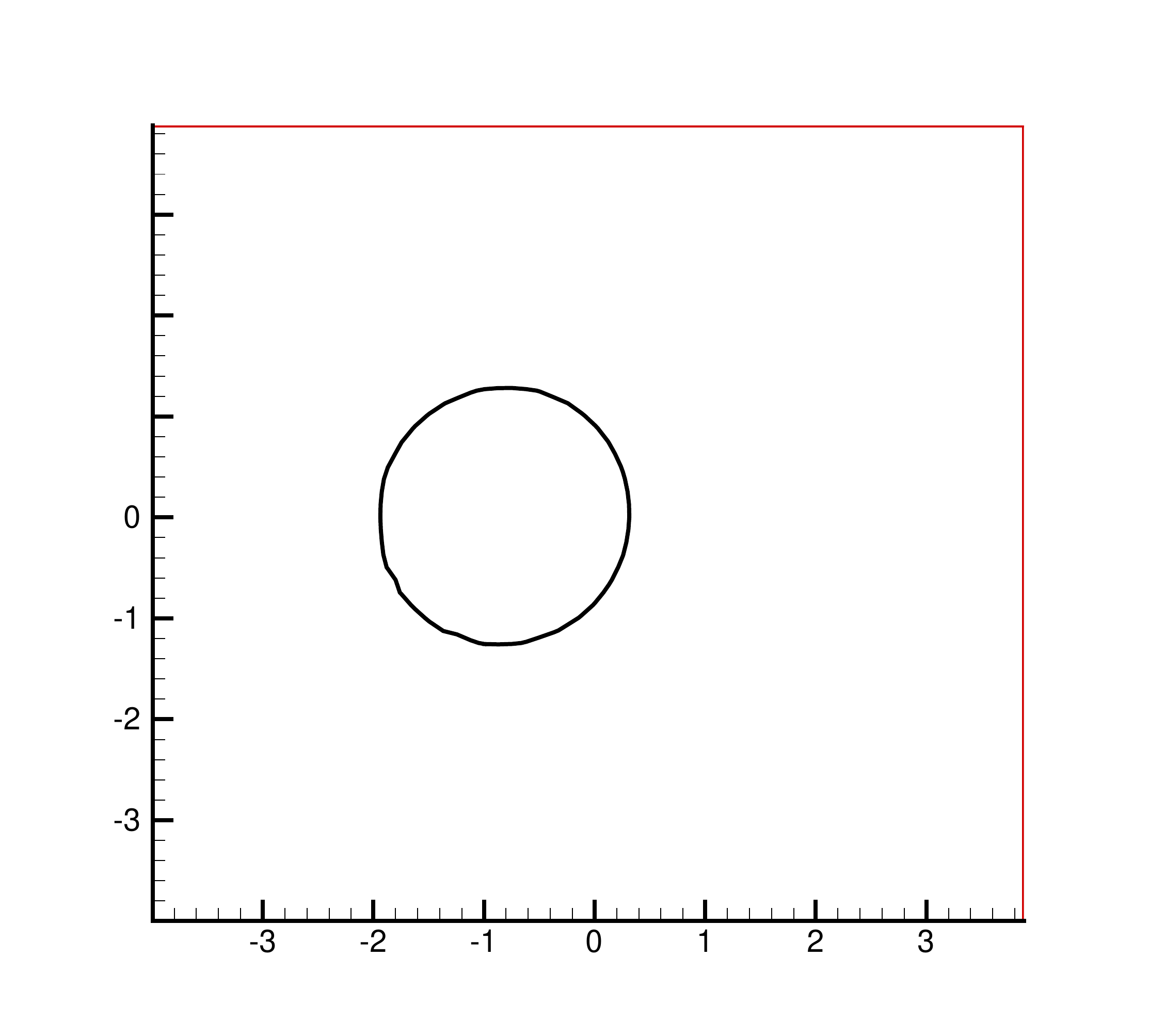}
}
\caption{Towards a circle.  {\it Evolution of Willmore of the
    regularized  flow with $\varepsilon=10h$ and $\varepsilon=10h^2$ when the initial level set
    corresponds to a singular curve.}}
\label{singular}
\end{figure}

{
We consider an initial data, where the initial level set is a singular curve.  We report the numerical results in
Figure \ref{singular},  which shows  the evolution of a singular curve into a circle. We approximate the Willmore flow with $\varepsilon=10h$ and $\varepsilon=10h^2$, respectively on a $64 \times 64$ grid for the computational domain $\Omega=[-4,4] \times [-4,4]$, which shows similar patterns with different regularization parameter $\varepsilon$.}

\subsection{Topology changes}

Finally, in this example, we present two examples to illustrate a topological change under Willmore flow, which is one of the major advantages of level set methods. However, level set methods for this problem suffer from an additional difficulty. The Willmore flow does not have a maximum principle, which leads to that two surfaces both undergoing an evolution by Willmore flow my intersect in finite time. To handle these difficulties, it is important to choose the regularization parameter $\varepsilon$ to avoid a blow up of the gradient of $\phi$ in finite time. In addition, it is important to choose small time steps to temporally resolve the rapid transients that occur at a topological change.

\begin{figure}[!ht]
 \centering
\subfigure[$t=0$] {
\includegraphics[width=0.3\columnwidth]{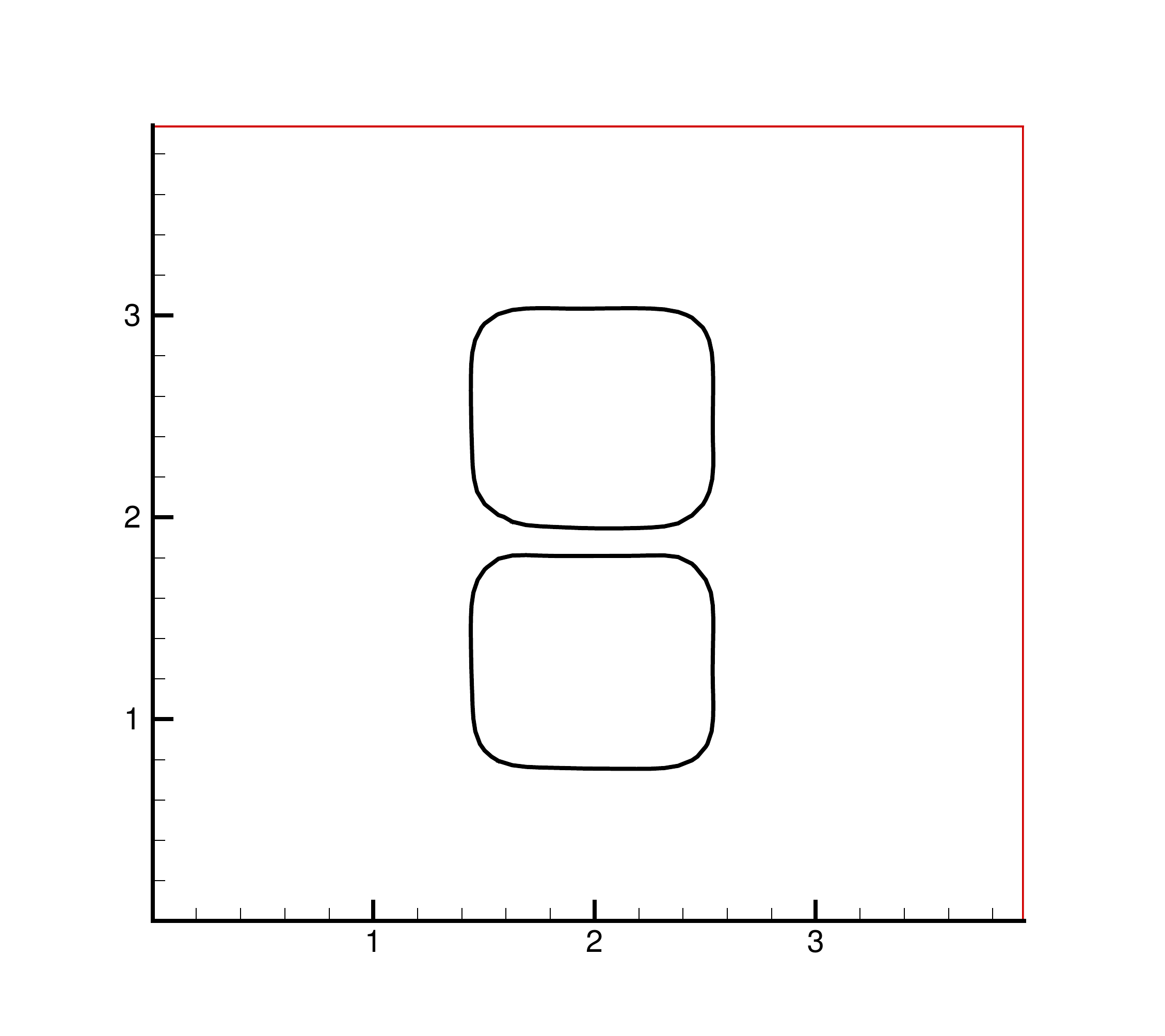}
}
\subfigure[$t=0.0002$] {
\includegraphics[width=0.3\columnwidth]{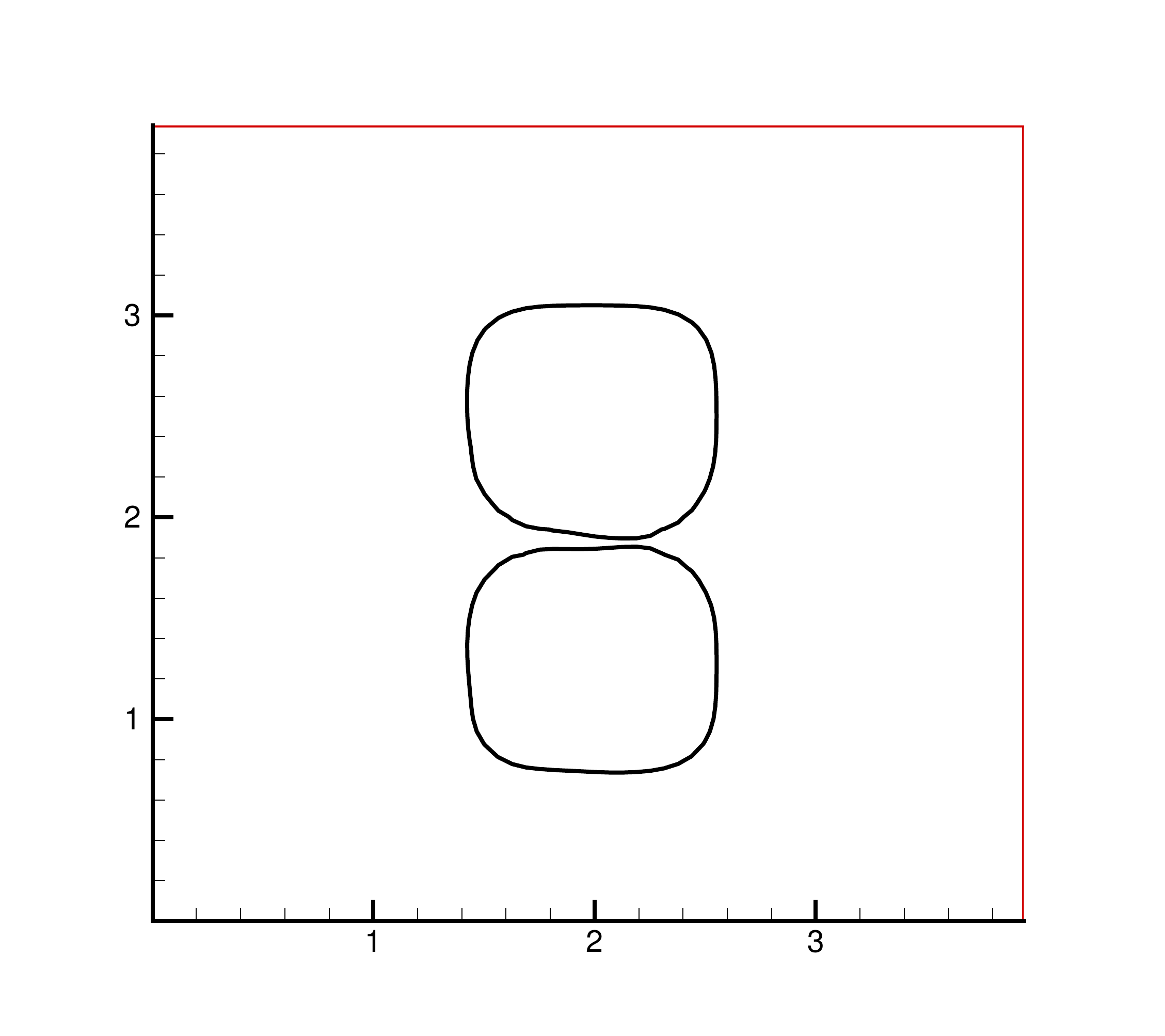}
}
\subfigure[$t=0.00022$] {
\includegraphics[width=0.3\columnwidth]{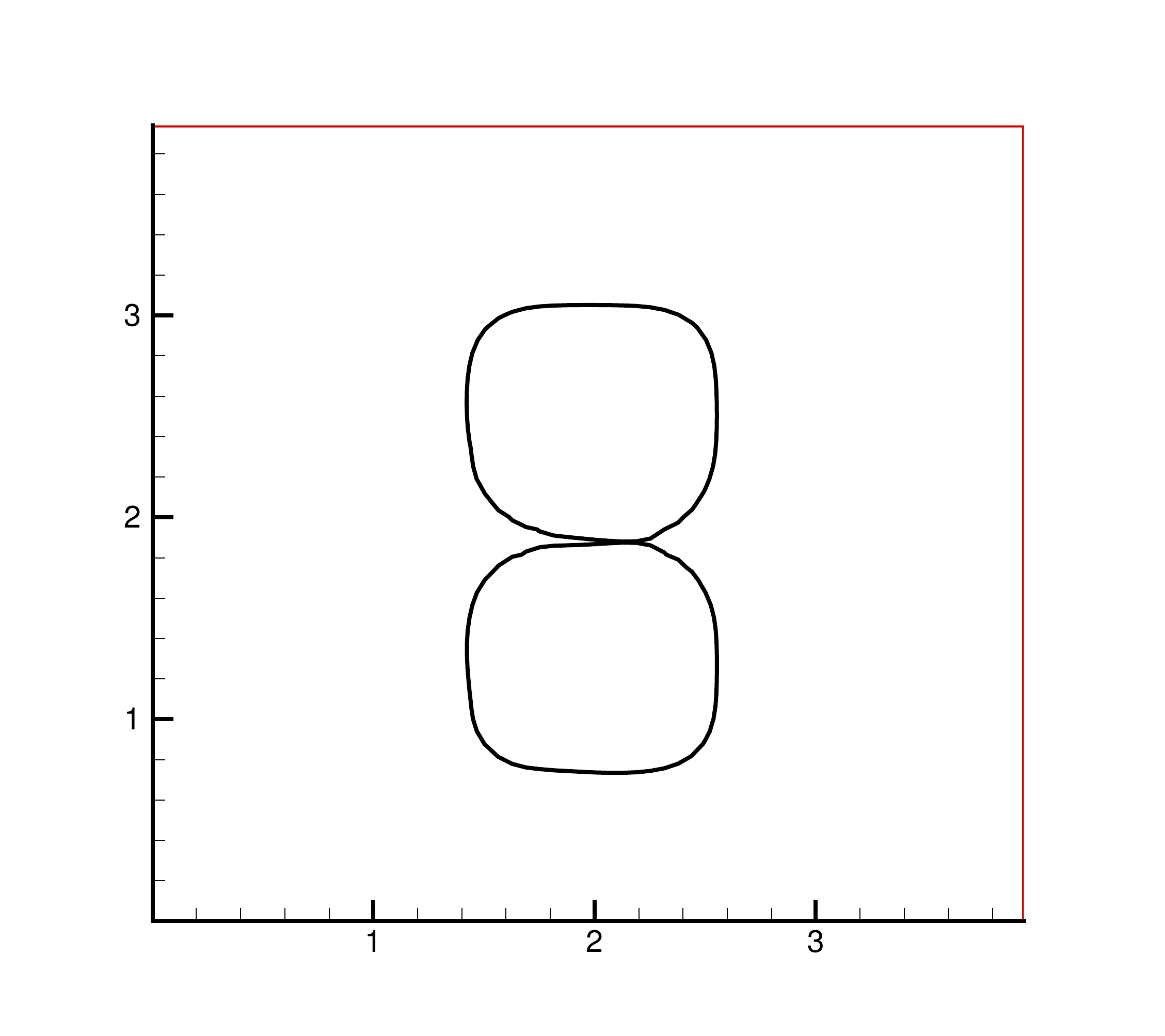}
}
\subfigure[$t=0.00023$] {
\includegraphics[width=0.3\columnwidth]{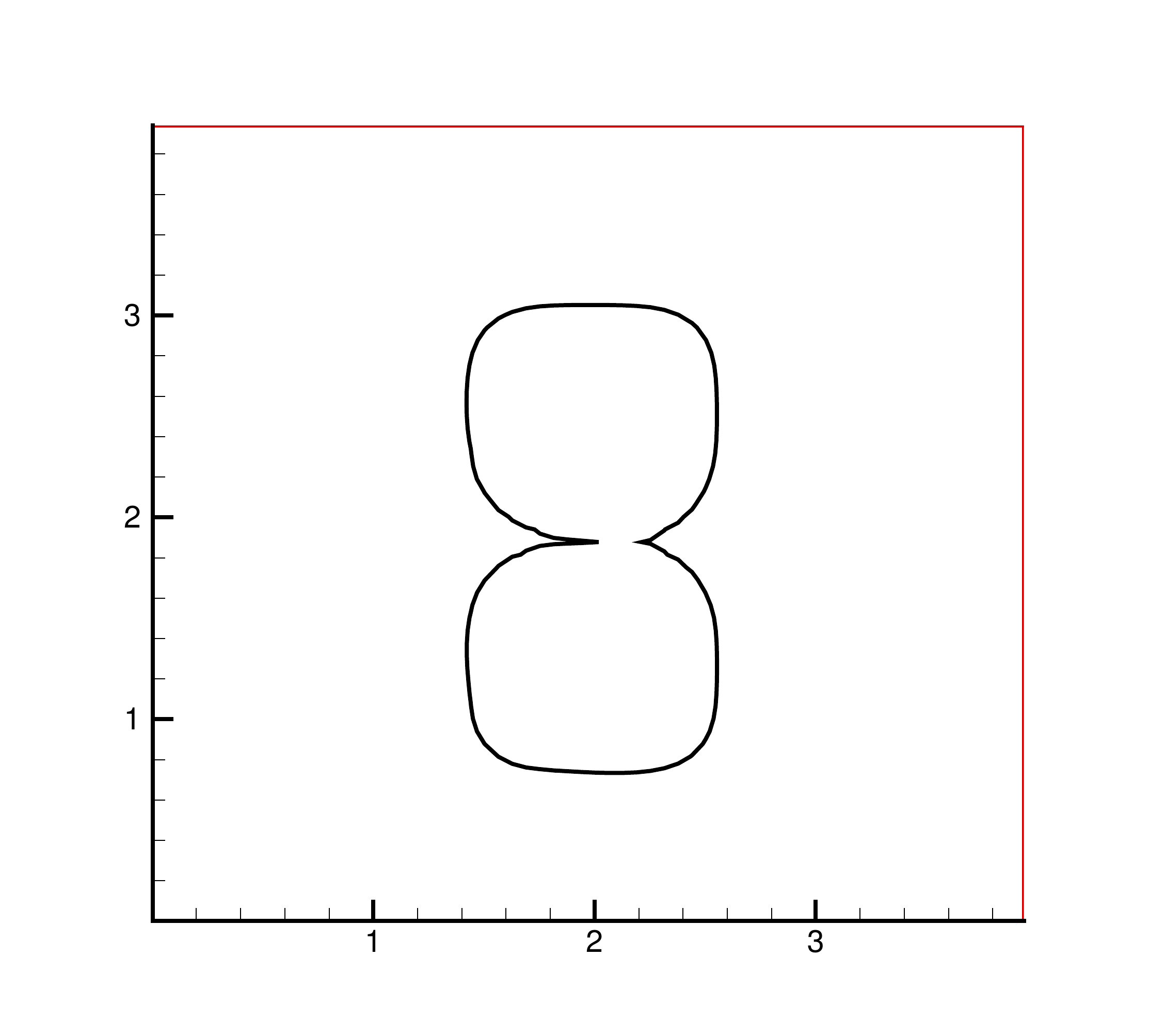}
}
\subfigure[$t=0.00025$] {
\includegraphics[width=0.3\columnwidth]{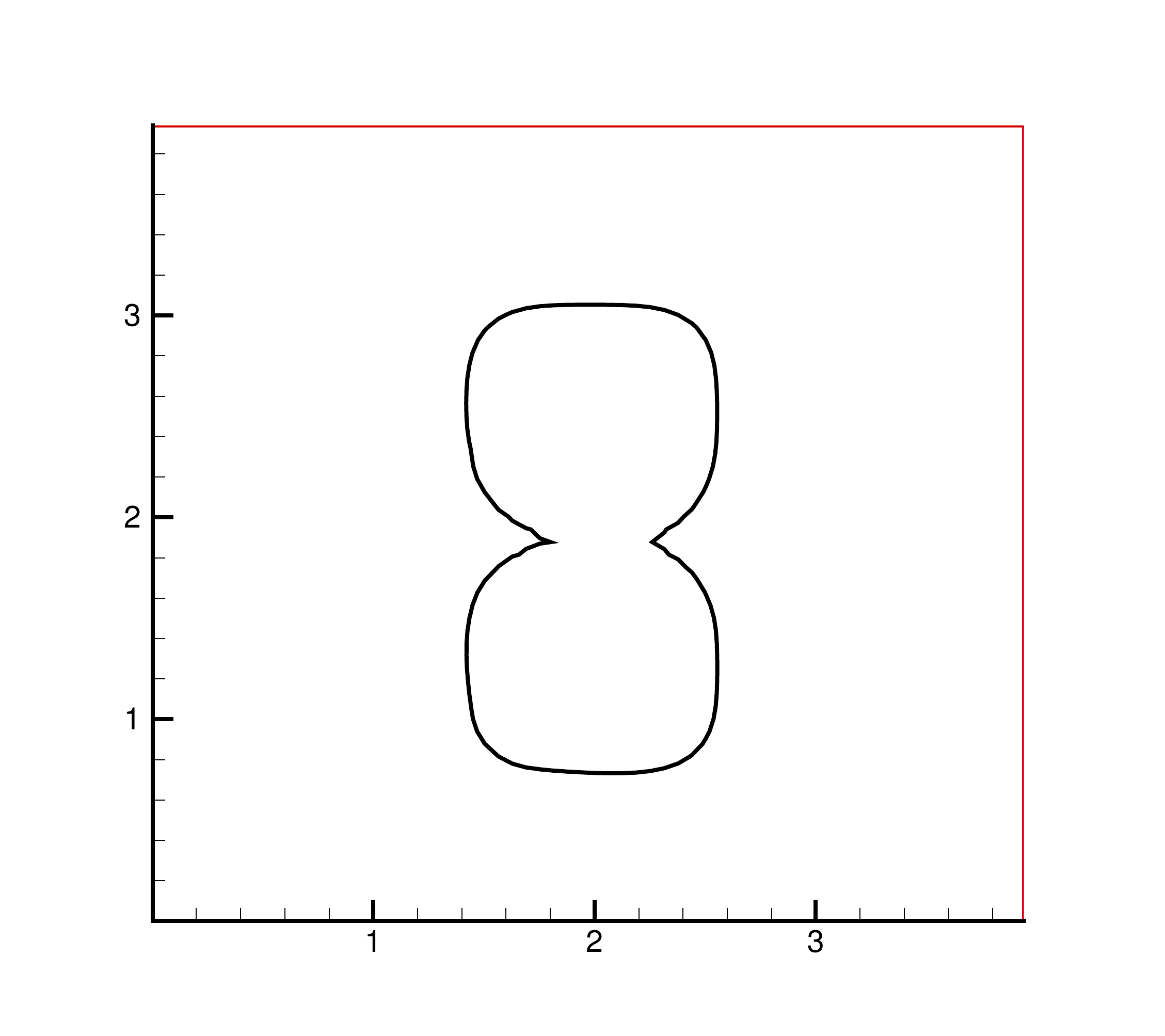}
}
\subfigure[$t=0.0003$] {
\includegraphics[width=0.3\columnwidth]{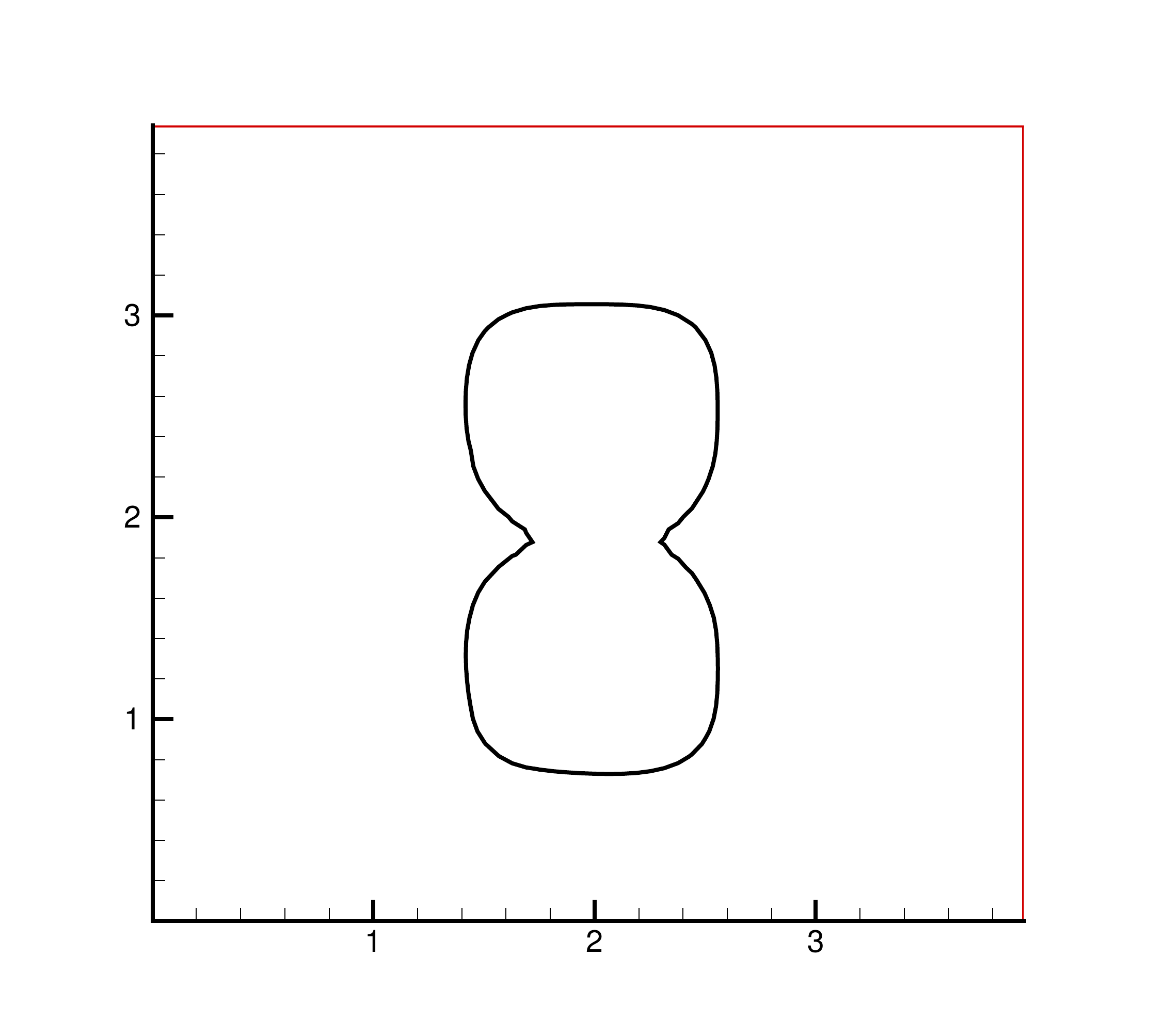}
}
\subfigure[$t=0.0007$] {
\includegraphics[width=0.3\columnwidth]{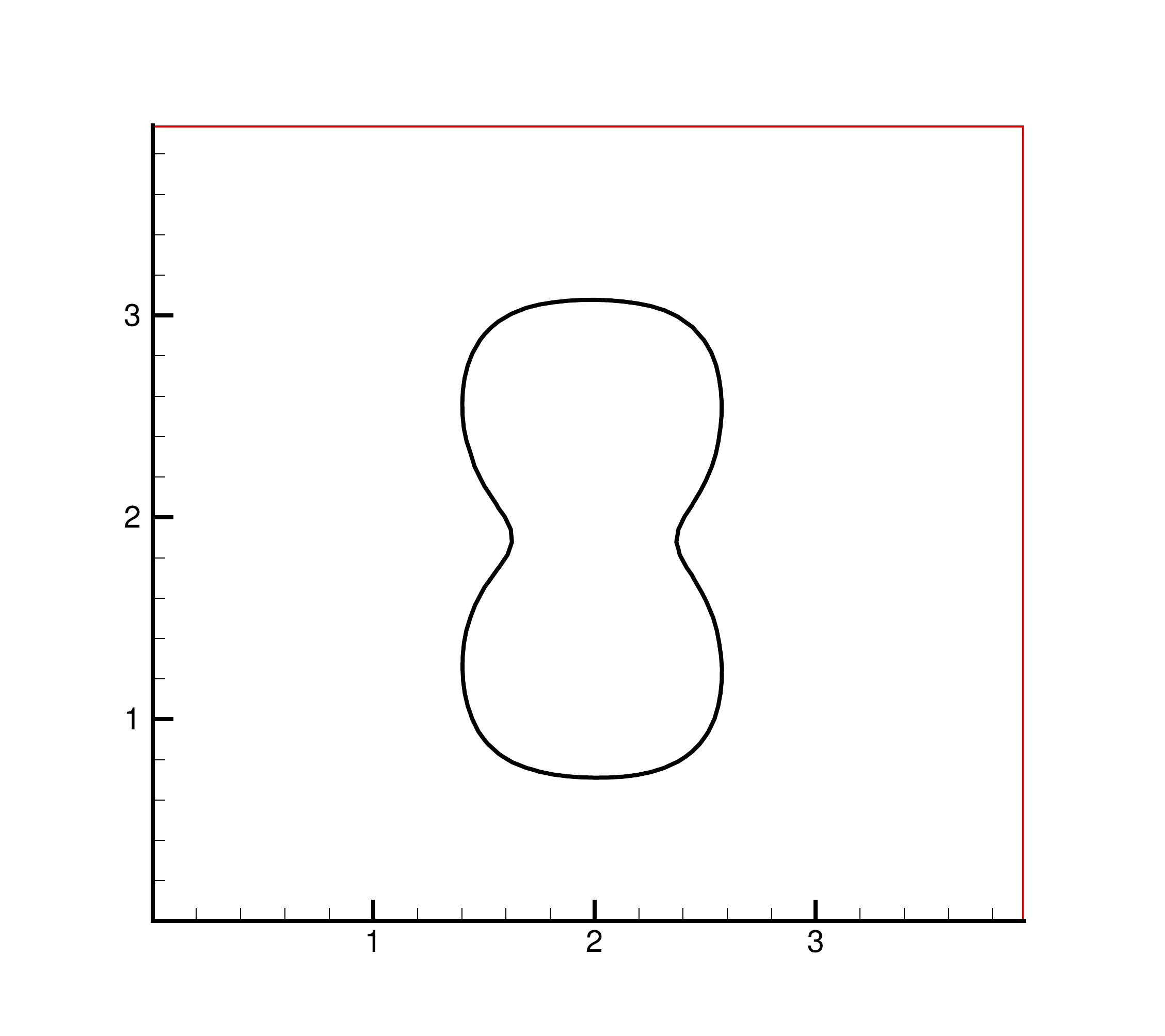}
}
\subfigure[$t=0.01$] {
\includegraphics[width=0.3\columnwidth]{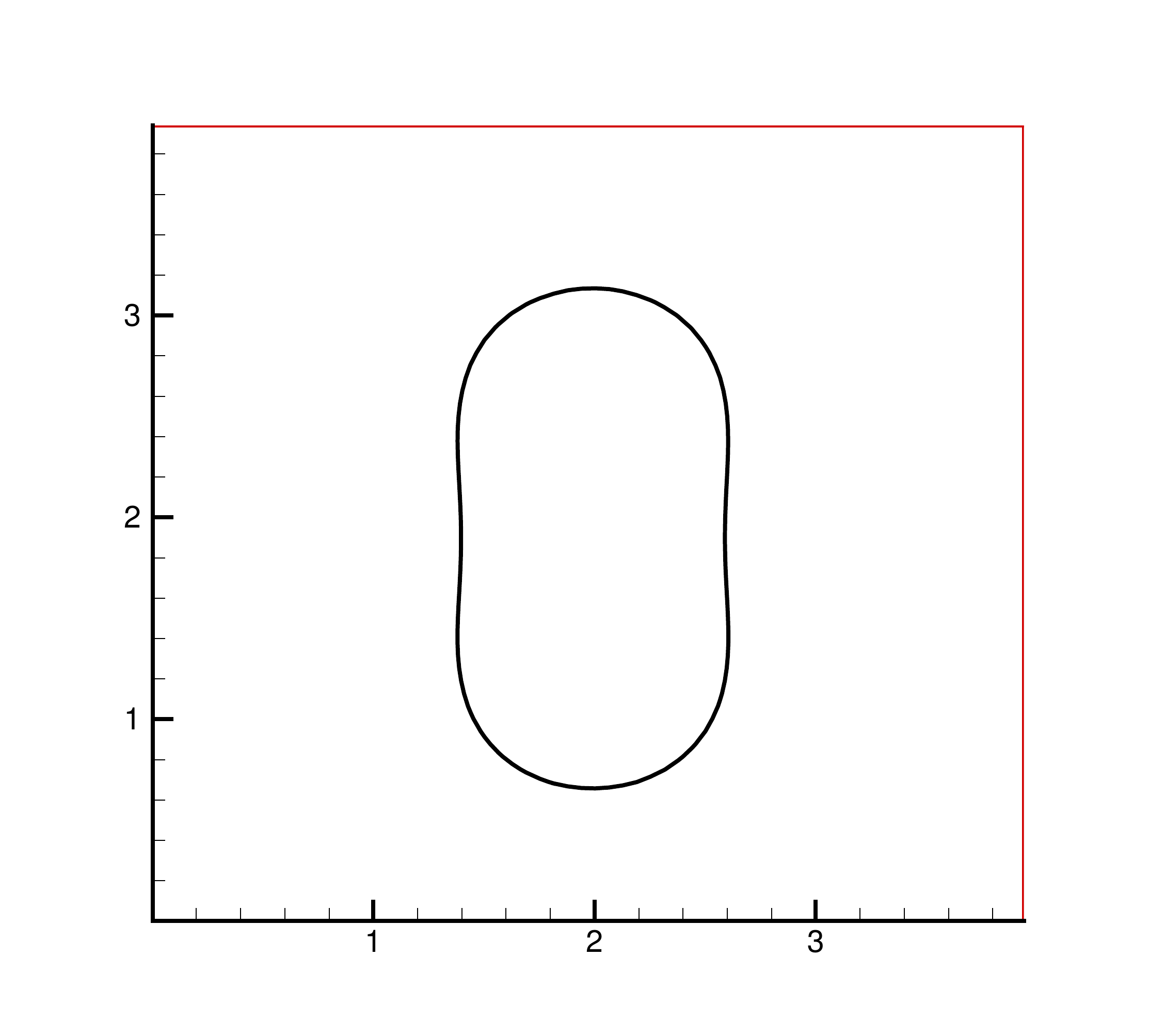}
}
\subfigure[$T=0.05$] {
\includegraphics[width=0.3\columnwidth]{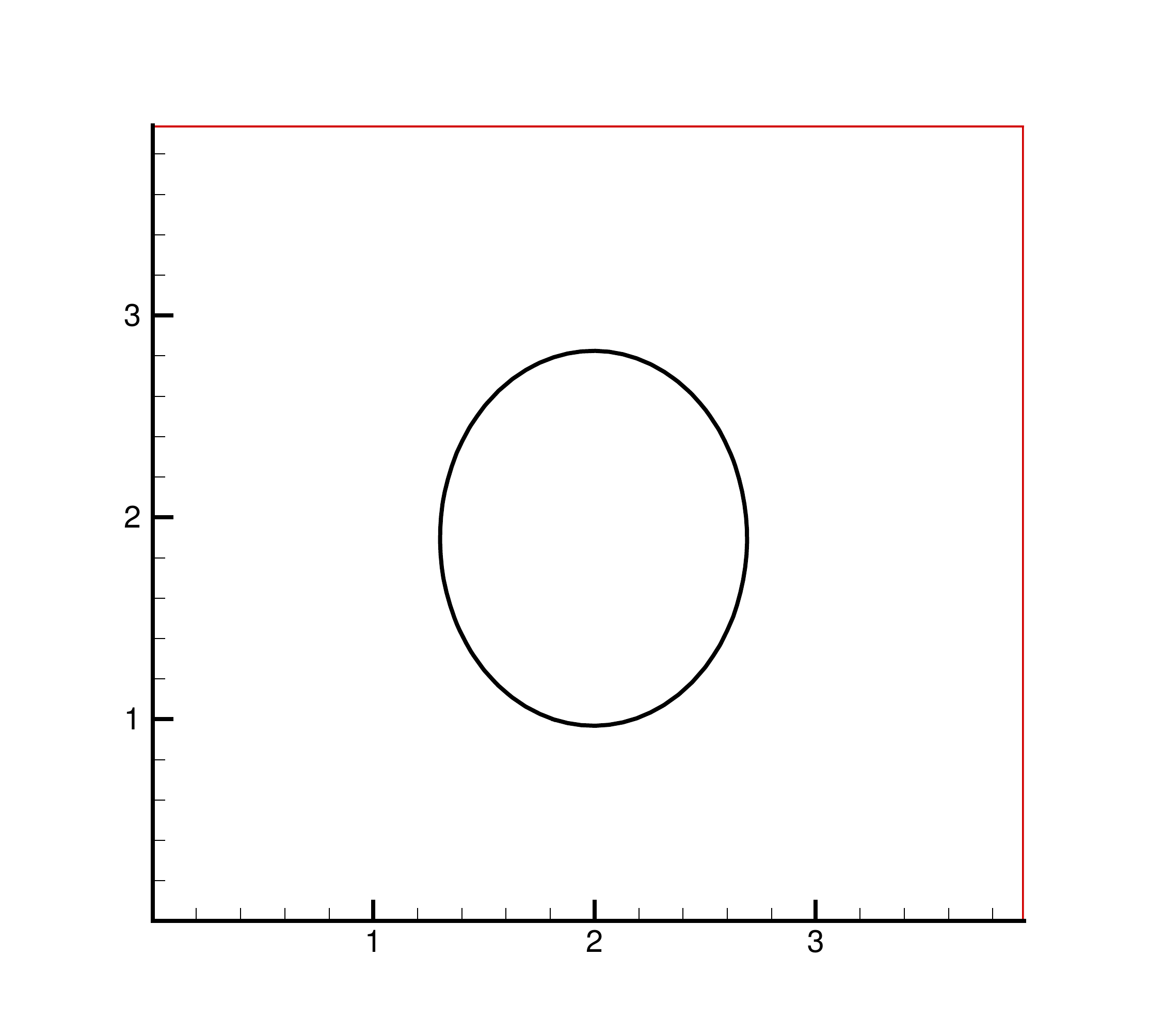}
}
\caption{Topology changes. {\it Two square-like shapes merge under the level set evolution of Willmore flow.}}
\label{two_square}
\end{figure}

Figure \ref{two_square} presents the merging of two square-like shapes
under the level set evolution of Willmore flow. We can see the two
square-like shapes merge and  eventually evolve to a single
circle. Here, we choose $\varepsilon=5h$, and the time step size is $\Delta t = 10^{-6}$. Since for such small time step, we use here a first order semi-implicit time marching method. { Certainly, we can use explicit Runge-Kutta method here, but because of the severe time step restriction of explicit methods, the time step must be very small, which takes about $\Delta t\approx 10^{-9}$.}

\begin{figure}[!ht]
\centering
\subfigure[$t=0$] {
\includegraphics[width=0.4\columnwidth]{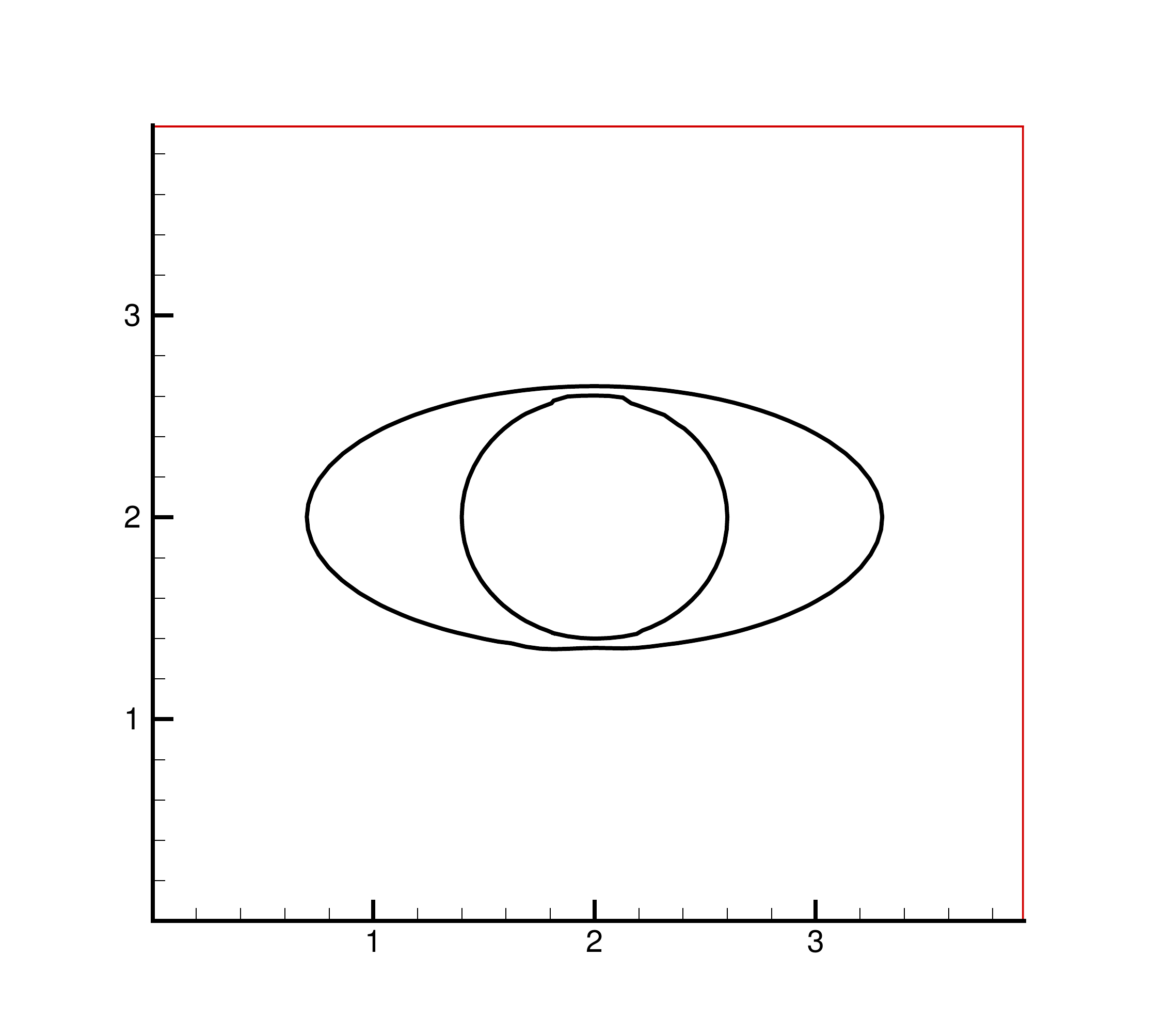}
}
\subfigure[$t=0.00002$] {
\includegraphics[width=0.4\columnwidth]{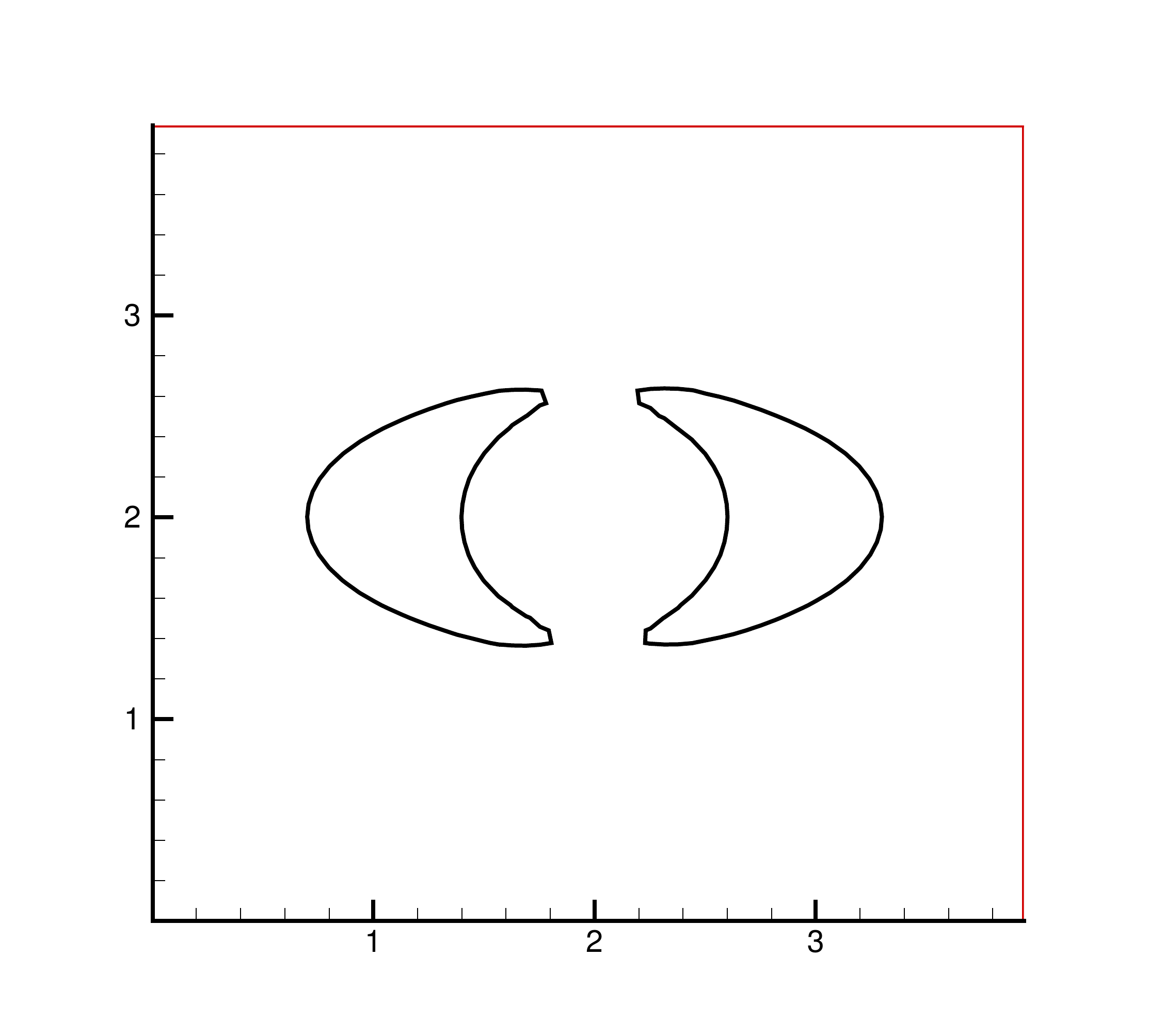}
}
\subfigure[$t=0.0001$] {
\includegraphics[width=0.4\columnwidth]{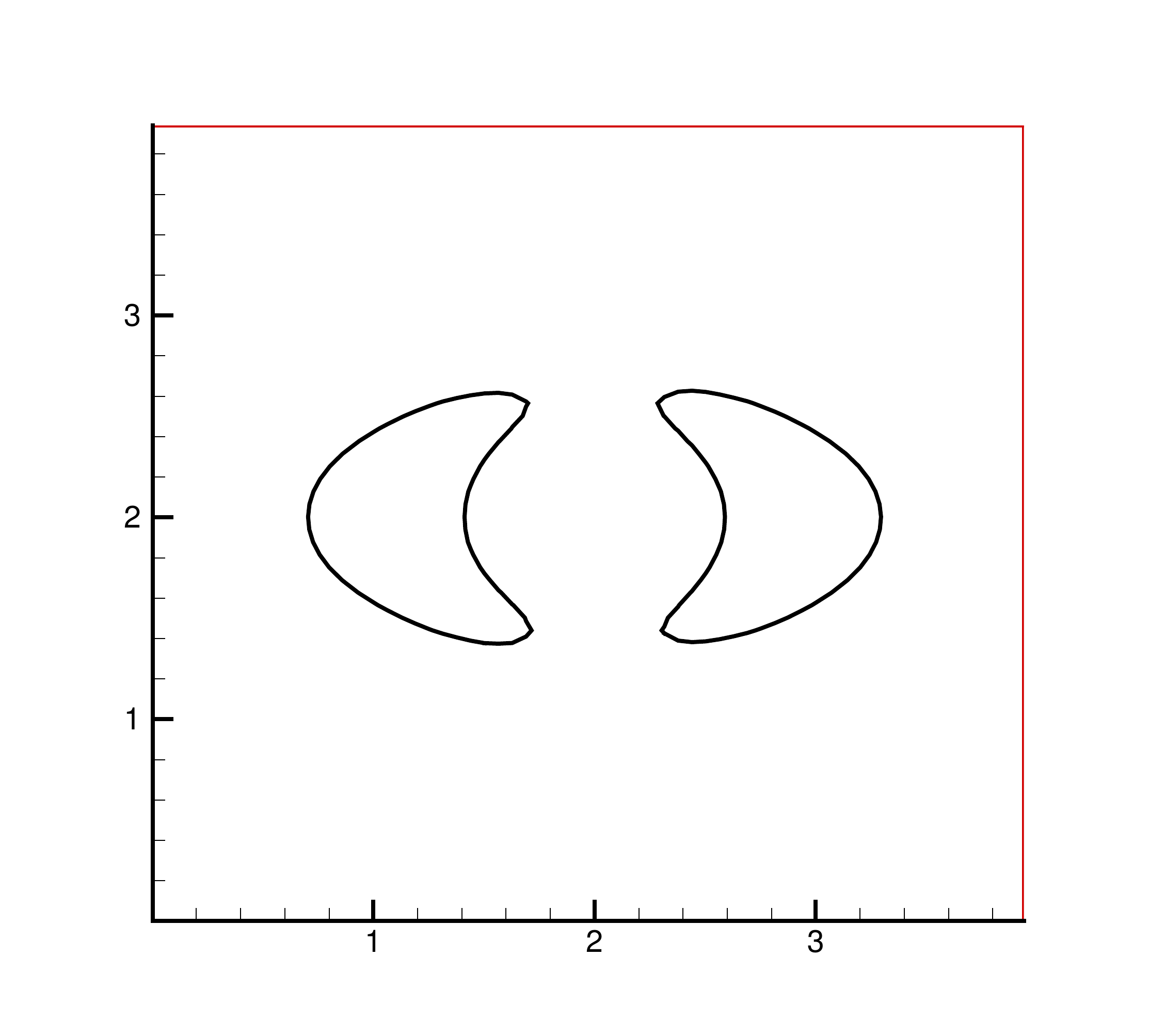}
}
\subfigure[$t=0.001$] {
\includegraphics[width=0.4\columnwidth]{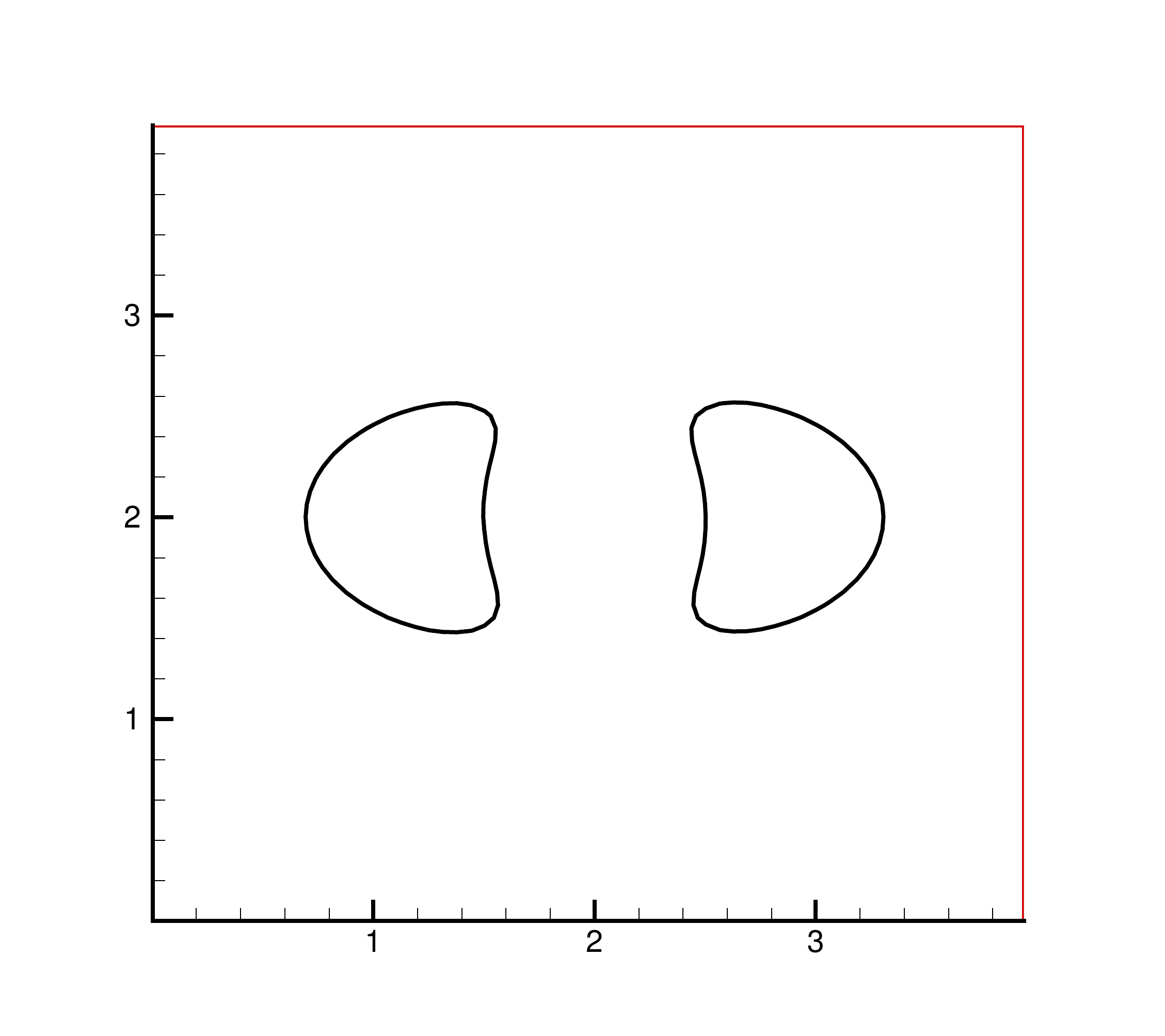}
}
\caption{Topology changes. {\it Two shapes pinch off under the level set evolution of Willmore flow.}}
\label{circle_ellipse}
\end{figure}

In Figure \ref{circle_ellipse}, we show another example of topological change. The initial configuration consists of two almost touching surfaces-the inner surface being a circle and the outer surface being an ellipse. We can see a phenomenon of pinching off in Figure \ref{circle_ellipse}. Here, we also choose $\varepsilon=5h$, and the time step size $\Delta t=10^{-6}$.

\section{Concluding remarks}
\label{se:con}
In this work, we have developed a local discontinuous Galerkin method for a level set formulation of Willmore flow.  The scheme is mass conservative and energy stable. In the level set method, we are mainly interested in capturing the movement of the interface. Therefore, it is reasonable and natural to combine our local discontinuous Galerkin scheme with $p$-adaptive technique to improve the efficiency of the method. To relax the severe time step restriction of explicit time integration methods for stability and achieve high order temporal accuracy, a semi-implicit Runge-Kutta method is adopted. The equations at the implicit time level are linear and we employ an efficient multigrid solver to solve them. Numerical experiments demonstrate that the proposed methods are efficient for Willmore flow, and the level set method can easily handle topological changes.

\begin{flushleft} \signFF \end{flushleft}
\begin{flushright} \signRG \end{flushright}

\end{document}